\documentclass[11pt, fleqn, a4paper]{article}
\usepackage{qtree}
\usepackage[utf8]{inputenc}
\usepackage{latexsym}
\usepackage{amsfonts}
\usepackage{amssymb}
\newtheorem{theorem}{Theorem}[section]

\newtheorem{settask}[theorem]{Task}

\newcommand{\nats}{{\mathbb{N}}}
\newcommand{\rats}{{\mathbb{Q}}}

\newcommand{\task}{\begin{settask}}
\newcommand{\etask}{\end{settask}}
\newcommand{\seqs}{{\mathbb{S}}}

\newcommand{\tail}[1]{#1'}

\newcommand{\newton}[1]{\underline{#1}}
\newcommand{\conv}[1]{#1^{\circ}}
\newcommand{\shi}[2]{#1^{-#2_\otimes}}
\newcommand{\shp}[2]{#1^{#2_\otimes}}
\newcommand{\co}{\circ}
\newcommand{\sbullet}{\mbox{\tiny \,$\bullet$\,}}
\newcommand{\rep}[1]{#1^{\tiny \sbullet}}

\newcommand{\set}{\mbox{\sf set\,}}

\newcommand{\lst}{\mbox{\sf list\,}}

\newcommand{\cyc}{\mbox{\sf cycle\,}}

\newcommand{\sfn}[1]{\mbox{\sf #1}}
\newcommand{\pdeg}{\mbox{\rm deg\ }}
\newcommand{\lne}{\mbox{\rm lgn}}
\newcommand{\gd}{\mbox{\rm gd}}
\newcommand{\hsp}{\hspace{3mm}}
\newcommand{\half}{\frac{1}{2}}
\newcommand{\binom}[2]{\left( \begin{array}{c}\hspace{-1ex} #1 \\\hspace{-1ex} #2 \end{array} \hspace{-1ex}\right)}
\newcommand{\perms}[2]{\left[ \begin{array}{c}\hspace{-1ex} #1 \\\hspace{-1ex} #2 \end{array} \hspace{-1ex}\right]}
\newcommand{\parts}[2]{\left\{ \begin{array}{c}\hspace{-1ex} #1 \\\hspace{-1ex} #2 \end{array} \hspace{-1ex}\right\}}

\newcommand{\falling}[2]{#1^{\underline{#2}}}
\newcommand{\rising}[2]{#1^{\overline{#2}}}
\newcommand{\reason}[1]{\mbox{\{#1\}}}

\newcommand{\leo}{\Lambda}
\newcommand{\loe}{\Lambda^{-1}}

\newcommand{\col}[1]{#1}
\newcommand{\uvec}[1]{\underline{#1}}

\newcommand{\revp}[1]{\grave{#1\,}}
\newcommand{\lde}[1]{#1}

\begin{document}
\def\qtreeinithook{\let\tabular=\oldtabular
\let\endtabular=\endoldtabular
\renewcommand{\tabcolsep}{6pt} }
%\title[Sequence Algebra]{In Praise of Sequence (Co-)Algebra and its implementation in Haskell}
%\author[K Clenaghan]{KIERAN CLENAGHAN}
\title{In Praise of Sequence (Co-)Algebra and its implementation in Haskell}
\author{Kieran Clenaghan}
\date{}
\maketitle
\begin{abstract}
What is Sequence Algebra?  This is a question that any teacher or student of mathematics or computer science can engage with.  Sequences are in Calculus, Combinatorics, Statistics and Computation.  They are foundational, a step up from number arithmetic.  Sequence operations are easy to implement from scratch (in Haskell) and afford a wide variety of testing and experimentation.  When bits and pieces of sequence algebra are pulled together from the literature, there emerges a claim for status as a substantial pre-analysis topic. Here we set the stage by bringing together a variety of sequence algebra concepts for the first time in one paper.  This provides a novel economical overview, intended to invite a broad mathematical audience to cast an eye over the subject.  A complete, yet succinct, basic implementation of sequence operations is presented, ready to play with.  The implementation also serves as a benchmark for introducing Haskell by mathematical example.  
\end{abstract}

\section{Introduction}
Consider these titles: \textit{Formal Power Series} \cite{Niven69}, \textit{Power Series, Power Serious} \cite{McIlroy99}, \textit{A Coinductive Calculus of Streams} \cite{Rutten05}, and \textit{Concrete Stream Calculus} \cite{Hinze10}. The mixing of the classical and the modern in these papers is stimulating and suggests a re-telling of the elementary theory and application of sequences.  Casting our net wider than the citations in those four papers, brings up a number of corroborating works, including \cite{Ward36,Brand64,Traub65,Tutte75}, in which the authors call attention to the intrinsic qualities and utility of an elementary calculus or algebra of sequences.  We do more than re-advertise this work -- we endeavour to tease out the common ground, emphasising economy of statement and notation, whilst embracing variety of approach.
     
Our aim is to attract those who are less well acquainted with sequence work, or those who are unfamiliar with Haskell or both.  It is so that others can enjoy ``messing about'', as Hayman \cite{Hayman91} might put it, with sequences and their implementation.   Elementary sequence algebra provides a good answer to the question, `What is the smallest coherent chunk of mathematics to set undergraduates to implement, from scratch, so that they get the greatest reward?'.

First we must say that sequence algebra is an umbrella title for algebraic manipulations of finite and infinite sequences, $p=[p_0,p_1, \ldots, p_n],\ f=[f_0,f_1,\ldots]$, over some given element set, $F$.  It encompasses, prominently, formal power series algebra, but is not restricted to it.  A finite sequence, viewed as a sequence of coefficients for powers of $x$, can be \textit{expressed} as a formal polynomial.  Thus $[0,1]=0x^0+1x^1=x$, or we may say that $x$ is \textit{implemented} by $[0,1]$.  Similarly, $3x^2-x+4=[4,-1,3]$.  A finite sequence can also be interpreted as an infinite sequence by appending an infinite sequence of zeros.  An infinite sequence can be \textit{expressed} as a formal power series:
\[ f=\sum_{i=0}^{\infty} f_ix^i = \sum_i f_ix^i= [f_0,f_1,f_2, \ldots]\]
We write $f_n$ or $[x^n]f$ for the element at position $n$, called the $n$th term of $f$; it is the coefficient of $x^n$ in the power series view.  The zeroth term can also be written $f(0)$, but otherwise the notation $f(n)$ is reserved for function application: let $p=3x^2+4$, then $p_2=3$, but $p(2)=16$; contrast $p_0=p(0)=4$.  This example reveals that the symbol $p$ is overloaded, it stands for a function of type $\nats \rightarrow F$ and also for another of type $F \rightarrow F$, and we must be careful to distinguish them.

It seems that in spite of Niven's paper \cite{Niven69} receiving the Lester R Ford award, the algebra of formal power series (FPS) has not entered the standard introductory university texts in any substantial way.  Niven's paper has, however, been a key reference for the first chapters of both \cite{Henrici74,GoJa83}, neither of which is standard undergraduate fare.  In general, we find that elements of FPS/sequence algebra appear throughout a vast literature, but the algebra itself is treated minimally, perfunctorily, or is taken for granted, as might be expected in research work \cite{EvPoShWa03} and specialist texts \cite{Niederhausen10,Basler00}.  It also comes under \textit{generatingfunctionlogy} \cite{Wilf90, GrKnPa94}.  The term ``generating function'' has proven to be a bit awkward, because of the ``function'' part.  Notice how many times Wilf \cite{Wilf90} and others have to remind readers when convergence is not as issue.  So, where some might use ``generating function'', we use the plain ``sequence expression''.  The term ``generating function'' is more applicable when an analytic function interpretation is intended \cite[Ch. VII]{Eil:74}.  However, we freely use the standard \textit{names} of core analytic functions for their Taylor sequences, but rather than say, $\exp(x)$, we use just $\exp$ for the sequence. 

The extent of material that fits into elementary sequence algebra is perhaps under-appreciated.  Our goal is to raise appreciation through a modest ``survey'' of examples, presented in section \ref{sec-survey}.  Various notations and proof-styles appear in the literature, and an effort has been made to be inclusive and to harmonise.  The one-word term ``sequence'' is often preferred to the three-word ``formal power series'', but both have their merits, the latter being preferred in the multivariate case, and when formal variable substitution is involved. 

The sequence algebra we exhibit is on the same elementary level as Niven's paper \cite{Niven69}.  Much detail has to be omitted so that we can cover more examples.  Enough detail is included to convey the foundational concreteness of the topic, and omitted detail is in the literature.    The Haskell implementation is given in full in section \ref{sec-Haskell} so that the reader can be in no doubt about how succinct it is, and can type it all up (or download it), and ``own'' it.  The more mathematically-inclined reader may dwell on section \ref{sec-survey}, and reflect on the potential of a sequence algebra topic in the mathematical curriculum.  The more programming-inclined may dwell on sections \ref{sec-Haskell}, \ref{sec-exercising} and \ref{sec-extensions}, and engage with the proposition that sequence algebra provides an excellent vehicle for exploring learning-mathematics-through-programming, or vice-versa.  There is great scope for making a contribution to the consolidation, refinement, and application of sequence algebra as an introductory subject.

To help catch, at a glance, some of the things that come under sequence algebra, we have included a number of tables.  Tables \ref{tab-diff-rules}, \ref{tab-integ-rules} and \ref{tab-logexp-props} are instantly recognisable as belonging to a calculus text, but here, uncommonly, the objects being related are sequences, not analytic functions.  Of course, they herald identities that hold for analytic functions, but in agreement with Niven \cite{Niven69} and Tutte \cite{Tutte75}, there is something to celebrate in the fact that the identities can be established on very elementary grounds.  There is more to celebrate in tables \ref{tab-seq-exps} and \ref{tab-biv-seqs}, because many of the sequence expressions therein have a dual interpretation as a set-theoretic structure specification \cite{FlSe09}.  Yet more satisfaction is to be had from table \ref{tab-core}, because the solutions to the defining differential equations for the core sequences transliterate trivially into Haskell definitions \cite{McIlroy01}.

Therefore a grounding in sequence algebra and its implementation surely pays off.  This claim is validated at least in the study of the two texts, \textit{Concrete Mathematics} \cite{GrKnPa94} and \textit{Analytic Combinatorics} \cite{FlSe09}.  There we see numerous examples where sequence algebra is at play, and the implementation can be used for testing, for reinforcement, or just for fun.  Some of these appear in the next section.  For example, in item {\bf W} we derive the bivariate power series $S(z,u)=\displaystyle \frac{\exp\co uz -1}{\exp\co z -1}$ that appears in \cite[sect. 7.6]{GrKnPa94}.  This generates a sequence $\breve{S}$ such that  $\breve{S}_m$ is a polynomial, and $\breve{S}_m(n)=\displaystyle\sum_{0\leq k<n}k^m$.  That is, $\breve{S}=\displaystyle[x, -\half x+\half x^2, \frac{1}{6}x-\half x^2 + \frac{1}{3}x^3, \ldots ]$.  It is striking how accessible the mathematics behind $S(z,u)$ is, and how easy it is to define the infinite $S(z,u)$ and $\breve{S}$ in a program.
\section{The basics}
Two sequences $f$ and $g$ are equal if they are equal at all indices: $\forall n\geq 0. [x^n]f=[x^n]g$.  There are many instances when it is obvious that a statement is subject to universal quantification, and in such cases we leave the $\forall$ part to be inferred by the reader.  Becoming fluent with the clumsy-looking \textit{coefficient extraction} operator, $[x^n]$, pays dividends \cite{GoJa83,Knuth94}. It obeys the precedence rules, $[x^n]f+g=([x^n]f)+g$, $[x^n]f g = [x^n](f g)$, and $[x^n]f\co g = [x^n](f\co g)$.  Moreover, it is  a linear operator:
\[ [x^n](f+g) = [x^n]f + [x^n]g; \hsp [x^n]cf = c[x^n]f \hsp \hsp \mbox{$c$ is a constant} \]  
The generalisation $[u^mz^n]$ will be used to identify a term in a bivariate sequence.  Let $E$ be the \textit{tail} or \textit{shift-left} operator: $E [f_0,f_1,f_2, \ldots]= [f_1, f_2, \ldots]$.  In formal power series language, $E\,f=\displaystyle\frac{1}{x}(f-f_0)$, and $f=f_0+xE\,f$; this is referred to as the \textit{head-tail property} (in \cite{Rutten05} it is the ``fundamental theorem of stream calculus'', see section \ref{sec-survey}, item {\bf P}).  We also have $[x^{n+1}]f=[x^n]E\,f$ and the head-tail expansion rule, $[x^n]f = [x^0]E^n f$.  We are free to mix notations -- here is the definition of convolution product, $f*g$, in which, typically, the $*$ is suppressed:
\[ [x^n]f g = \sum_{k=0}^n f_k g_{n-k} \]
Observe that, since $x_1=1$ is the only non-zero element of $x=[0,1]$, we have
\[ [x^0]x f = 0; \hsp [x^n]x f = x_1f_{n-1} = f_{n-1}, \hsp \mbox{for\ } n>0 \]
Thus $x[f_0,f_1,f_2,\ldots ] = [0,1]*[f_0,f_1,f_2,\ldots ] = [0,f_0,f_1,f_2, \ldots]$, and $x$ can be viewed as a \textit{right-shift} operator.  The absorption law, $[x^n]x^mf = [x^{n-m}]f$, is simple but effective.  The product in both $f=f_0+xE\,f$ and $f= \sum_i f_ix^i$ is convolution product, $+$ is pointwise addition, the element $f_i$ is automatically identified with the singleton sequence, $[f_i]$, as required, and $x=[0,1]$.  A recursive equation for product is easily derived (and just as easily translated into Haskell); $E\,f$ is abbreviated to $f'$:
\begin{equation}
f g =  (f_0 + xf')(g_0 + xg') =  f_0g_0 + f_0xg' + xf'g  = f_0g_0 + x(f_0g'+f'g) \label{eq-conv-prod}
\end{equation}
Sequences, $\seqs$, over an integral domain or field $F$ (known from context) form an integral domain \cite{Niven69,Lipson81} ($F[[x]]$ is the standard notation for formal power series in $x$ over $F$).  We will keep $F=\rats$ in mind.  The subsets $\seqs_0,\ \seqs_1,$ and $\seqs_{\neq\! 0}$ comprise sequences with zeroth term 0, 1, and non-zero, respectively.  A subset $\seqs_C$ of $\seqs_0$ comprises sequences $f$ in which $f_1\neq 0$, that is $E\,f \in \seqs_{\neq\! 0}$.  Unique square (and $n$th) roots exist for sequences in $\seqs_1$.  Sequence composition $f\co g = \sum_k f_kg^k$ is defined for $g\in \seqs_0$. A unique compositional inverse, $\conv{f}$, called \textit{converse}, exists for $f\in \seqs_C$.  The notation follows {\cite{BiDM97} and distinguishes converse $f\co \conv{f}=x$ from multiplicative inverse $f*f^{-1}=1$.  The latter exists for $f\in \seqs_{\neq\! 0}$.  Differentiation, $D$, is term-wise, as for formal polynomials: $[x^{n}]D\,f = (n+1)[x^{n+1}]f$.  A little induction gives the Maclaurin expansion rule: $[x^n]f=\displaystyle\frac{1}{n!}[x^0]D^nf$.  From the definition of $\int$ as a \textit{right} inverse to $D$, $[x^n]D \int f  = [x^n]f$, we deduce $[x^{n+1}]\int f = \displaystyle\frac{1}{n+1}[x^n]f$.  Setting $[x^0]\int f  =  0$, the fundamental theorem of sequence calculus (FTC) is immediate:
\[ f   =  f_0+\int D f \hspace{10mm} D(\int f)   =  f \]
\begin{table}
\[ \begin{array}{lrcll}
1.& D\,a & = & 0 & \mbox{\bf constant}\\
2.& D\,x & = & 1 & \mbox{\bf variable}\\
3.& D\,(f+g) & = & D\,f+D\,g & \mbox{\bf sum}\\
4.& D\,(f g) & = & (D\,f) g + f (D\,g)& \mbox{\bf product}\\
5.& D\,(f^{-1}) & = & (-D\,f)/f^2 & \mbox{\bf reciprocal}\\
6.& D\,(f/g) & = & ((D\,f)g-f(D\,g))/g^2 & \mbox{\bf quotient}\\
7.& D\,(\conv{f}) & = & ((D\,f) \co \conv{f})^{-1} & \mbox{\bf converse}\\
8.& D\,(f^n) & = & n(D\,f)f^{n-1} & \mbox{\bf power}\\
9.& D\,(a_nx^n) & = & na_nx^{n-1} & \mbox{\bf monomial}\\
10.& [x^n]D & = & (n+1)[x^{n+1}] & \mbox{\bf power-series}\\
11.& D\,(f \co g) & = & ((D\,f) \co g)(D\,g) & \mbox{\bf composition}\\
12. & [x^n]f & = & \displaystyle\frac{1}{n!}[x^0]D^nf & \mbox{\bf Maclaurin}
\end{array}\]\caption{differentiation rules\label{tab-diff-rules}}
\end{table}
The familiar rules in tables \ref{tab-diff-rules} and  \ref{tab-integ-rules} have easy sequence-algebraic proofs (the third $\int$-product rule is called the \textit{differential Baxter axiom} in \cite{BuRo12}).  For example, here is a typical proof of the differential composition rule, or chain rule \cite[Ch. 1]{Henrici74}.  A proof by coinduction is presented for contrast in item {\bf Q} of the next section.  Let $f\in \seqs,\ g\in \seqs_0$, then,
\begin{eqnarray*}
\ D\,f\co g & = & f_1+2f_2g+3f_3g^2+\cdots +kf_kg^{k-1}+\cdots \\
\ (D\,f\co g)D\,g & = & f_1D\,g+2f_2gD\,g+3f_3g^2D\,g+\cdots +kf_kg^{k-1}D\,g+\cdots \\
  & = & \reason{the power rule, $kg^{k-1}D\,g=D\,g^k$} \\
  & & f_1D\,g + f_2D\,g^2 + f_3D\,g^3+ \cdots \\
\  [x^n](D\,f\co g)D\,g &=& \reason{distribute $[x^n]$, use $[x^n]D\,g^k=(n+1)[x^{n+1}]g^k$}\\
&& (n+1)[x^{n+1}](f_1g^1+f_2g^2+f_3g^3+\cdots+f_{n+1}g^{n+1} ) \\
  & = & (n+1)[x^{n+1}]f\co g \\
   & = & [x^n]D(f\co g)
\end{eqnarray*}  
\begin{table}
\[ \begin{array}{lrcll}
1.&\displaystyle \int a_nx^n & = &\displaystyle \frac{a_n}{n+1}x^{n+1} & \mbox{\bf monomial}\\
2.&\displaystyle\int(af + bg) & = & \displaystyle a\int f + b\int g & \mbox{\bf linear} \\
3.&\displaystyle(n+1)[x^{n+1}]\int & = & \displaystyle[x^n]  & \mbox{\bf power series} \\
4.&\displaystyle\int ((D\,f) \co g)aD\,g & = &\displaystyle a f\co g + c &\mbox{\bf composition} \\
5.&\displaystyle\int f D\,g & =&\displaystyle fg - \int (D\,f) g + c &\mbox{\bf product (1)} \\
6.&\displaystyle\int fg & = & \displaystyle f \int g - \int ((D\,f) \int g)& \mbox{\bf product (2)} \\
7.&\displaystyle\int (D\,f) \int (D\,h) + \int D\,(fh) & = &\displaystyle h\int D\,f +  f\int D\,h&\mbox{\bf Baxter}
\end{array}\]\caption{integration rules\label{tab-integ-rules}}
\end{table}
\section{Sequence Algebra examples}\label{sec-survey}
The following itemization ({\bf A-Z}) of snippets provides a brief survey of the character of sequence algebra.  It is, of course, only a small fraction of the subject.

\ \\
\noindent {\bf (A)} Sequence algebra is foundational in the sense that it is a low-level concrete extension of arithmetic.  To appreciate this, try making the following simpler.  Define $\exp$ by the sequence differential equation $D\,\exp = \exp;\ \exp_0=1$.  Then, by the Maclaurin rule, $\exp=1+x/1!+x^2/2!+x^3/3!+\cdots=[1,1,1/2,1/3!, \ldots ]$.  Define $x^*$ by the sequence difference equation $E\,x^*=x^*;\ x^*_0=1$;  then, by the head-tail property, $x^*=1+x+x^2+\cdots =[1,1,1,\ldots]=1/(1-x)$; the notation is based on the Kleene star \cite{Eil:74}.  Let $\log f=\lne\co (f-1)$ where $\lne$ (for which there is no established name other than $\log(1+x)$) is the converse of $\exp-1$; that is, $\lne=\conv{(\exp-1)}$.
Observe that $(\exp-1)\co \lne=x$ implies $\exp \co \lne = 1+x$, so
$ D\,\lne=D\,\conv{(\exp-1)} = (\exp \co \lne)^{-1}=1/(1+x)$. The FTC gives $\lne = \int 1/(1+x)$, and the coefficients can be calculated:
\begin{eqnarray*}
[x^{n+1}]\lne &=& [x^{n+1}]\int \frac{1}{1+x} = \frac{1}{n+1}[x^n]\frac{1}{1+x} = \frac{1}{n+1} (-1)^n\\  
\lne &=& x-\frac{x^2}{2}+\frac{x^3}{3}-\frac{x^4}{4}+\cdots
\end{eqnarray*} 
Some well-known rules relating to $\log$ and $\exp$, derivable from the differential equation for $\exp$ using sequence algebra,  appear in table \ref{tab-logexp-props} (preconditions are omitted to avoid clutter).  Most of these are meticulously proven by Niven \cite{Niven69}.  However, Niven does not use composition $(\co)$, but by using composition and its associativity and distributivity laws \cite{Henrici74}, his theorem 17 and proof can be rendered as follows: let $g=1+f,\ f,h\in \seqs_0$,
\begin{eqnarray*}
g^r\co h & = & \exp \co (r\log g) \co h = \exp \co \,r((\lne \co f) \co h) \\
 & = & \exp \co \,r (\lne \co (f \co h)) = \exp \co \,r \log ((f \co h)+1) \\
 & = & \exp \co \,r \log ((f+1) \co h) = \exp \co \,r \log (g\co h) = (g\co h)^r 
\end{eqnarray*}
Proof of Euler's identity, $\exp\co ix = \cos+i\sin$, illustrates appeal to the uniqueness of solution to certain differential equations.  Let $D\,\sin=\cos;\ \sin_0=0$, $D\,\cos =-\sin;\ \cos_0=1$, and $i^2=-1$.  Then both $\cos + i\sin$ and $\exp \co ix$ satisfy $D\,g=i\,g;\ g_0=1$, and therefore must be equal.  De Moivre's theorem follows:
\[ (\cos + i\sin)^n  =(\exp \co ix)^n  = \exp \co ix \co nx 
 = (\cos + i\sin) \co nx = \cos\co nx + i(\sin \co nx)\]
\begin{table}
\[\hspace{-10mm} \begin{array}{rclclll}
 \exp \co (\log f) &=& f &&& \mbox{\bf $\exp$ converse} \\
 \lne \co f & = &\lne \co g & \Rightarrow & f=g \hsp&\mbox{\bf $\lne$ cancellation} \\
 \exp \co f & = &\exp \co g & \Rightarrow & f=g &\mbox{\bf exp cancellation} \\
 \log f & = &\log g & \Rightarrow & f=g &\mbox{\bf log cancellation} \\
 D\,(\log g) & = & \displaystyle\frac{D\,g}{g} &&&\mbox{\bf log derivative} \\
 \lne\co f &=&0 & \Leftrightarrow & f=0 & \mbox{\bf zero $\lne$} \\
 \log g &=&0 & \Leftrightarrow & g=1 & \mbox{\bf zero log} \\
 \log (fg) & = & \log f + \log g &&&\mbox{\bf log product} \\
 \log f^r & = & r\log f &&& \mbox{\bf log rational power} &(r\in \rats) \\
 \exp\co(f+g) & = & (\exp\co f)(\exp \co g) &&& \mbox{\bf exp sum} \\
 \exp^n\co f & = & \exp \co nf &&& \mbox{\bf power exp} \\
 f^r & = & \exp \co (r \log f) &&& \mbox{\bf general power defn}&(r \in F)\\
 g^r\co h &=& (g\co h)^r &&&\mbox{\bf general distributivity}&(r \in F)\\
 f^rf^s & = & f^{r+s} &&& \mbox{\bf law of exponents}&(r,s \in F) \\ 
 D\, f^r & = & r f^{r-1}D\, f &&&\mbox{\bf differential $F$-power}&(r \in F)
\end{array}\]\caption{Rules relating to $\log$ and $\exp$ \label{tab-logexp-props}}
\end{table}
\noindent {\bf (B)} A \textit{counting} sequence, $c$, is either \textit{ordinary}, $c=[c_0,c_1,c_2,c_3\ldots]$ or \textit{exponential}, $c=[c_0,c_1/1!,c_2/2!, c_3/3!,\ldots]$.  In either case, $c_n$ counts the number of objects of size $n$ generated by some structure specification, $C$.  In the former case $c_n=[x^n]c$, and in the latter case $c_n=n![x^n]c$.  Roughly speaking, a structure is built from nodes, and its size is the number of nodes it has.  The nodes may be labelled or unlabelled.  Exponential sequences are used for labelled objects because, in that case, convolution product automatically counts all possible labellings in making an ordered product.  Let $f$ and $g$ count labelled structures (generated by some $F$ and $G$, respectively); then ordered pairs of such structures are counted by
\[ [x^n]fg = \sum_{k=0}^n \frac{f_kg_{n-k}}{k!(n-k)!}=\sum_{k=0}^n\binom{n}{k}f_{k}g_{n-k}/n! \]
The ordered $k$-fold product, $F^k$, has counting sequence $f^k$.  Ordered lists of $F$-objects are counted by 
$\lst \co f = x^* \co f = f^*$.  If the order does not count, then we use $\set \co f = \sum_{k} f^k/k! = \exp \co f.$
 
The fact that permutations can be written as sets of cycles can be made explicit in the definition of the counting sequence for permutations \cite{FlSe09,Cameron17}:
\[ \sfn{perm} = \sfn{set}\co\sfn{cycle} \]
Just put $\sfn{set}=\exp$ and $\sfn{cycle}=\log x^*$.  The sequence \sfn{perm} is $x^*$ regarded as an \textit{exponential sequence}, $x^*=[1,1!/1!,2!/2!,3!/3!\ldots]$.  Removal of the factorial divisors is performed by $\leo$, and $\leo \sfn{perm} = [1,1,2,6,24,120,\ldots]$, $[x^n]\leo \sfn{perm}=n![x^n]\sfn{perm}=n!$, the number of permutations of $n$ symbols.  There are $(n-1)!$ cyclic permutations of $n$ symbols, and the exponential counting sequence for these is:
\[ \cyc= \sum_{n>0} (n-1)!\frac{x^n}{n!} = \sum_{n>0}\frac{1}{n}x^n =  -\lne\co (-x)  =-\log (1-x)= \log  x^* \]     

\ \\
\noindent {\bf (C)} The number of ways, $s_n$, of inserting brackets into a list of $n$ symbols, subject to well-formedness, is counted by the Hipparchus-Schr\"{o}der sequence \cite{Stanley97,Stanley11}, $s=\displaystyle \frac{1}{4}(1+x-\sqrt{1-6x+x^2}\,)=[0,1,1,3,11,45,197,903,4279,20793,103049,\ldots]$.  This can be derived from the specification of a bivariate counting sequence for Schr\"{o}der trees:
\begin{equation}
 \sfn{schroeder(z,u)} = z+u*(\sfn{pluralList} \co \sfn{schroeder(z,u)}) \label{eq-schroder}
\end{equation}
Here $\sfn{pluralList} = \lst-x-1$.  The coefficient of $u^kz^n$ in $\sfn{schroeder(z,u)}$ gives the number of Schr\"{o}der bracketings of $n$ symbols using $k$ pairs of brackets.  Observe that equation (\ref{eq-schroder}) can also be read \cite{FlSe09} as a set-theoretic specification of Schr\"{o}der trees, with $z$ and $u$ naming different kinds of nodes.  Using $\sfn{pluralList}=x^*-x-1$, the above definition of $s$ derives, via the quadratic formula, from the equation $s=\sfn{schroeder(x,1)}=x+(x^*-x-1)\co s= x+s^*-s-1=\displaystyle x+\frac{1}{1-s}-s-1$.  Here is a foretaste of computing the Schr\"{o}der numbers in Haskell, which is explained in section \ref{sec-Haskell}:
\begin{verbatim}
    schroeder   =  z + u*(pluralList `o` schroeder)
    > takeBiv [1..6] schroeder
    [[0],[1,0],[0,1,0],[0,1,2,0],[0,1,5,5,0],[0,1,9,21,14,0]]
    > takeW 11 ((1+x-sqroot(1-6*x+x^2))/4)
    [0,1,1,3,11,45,197,903,4279,20793,103049]
\end{verbatim}   
This reveals, for example, $[u^2z^5]\sfn{schroeder}=9$, that is 9 bracketings of 5 symbols with 2 pairs of brackets.  One can see that the elements of the second sequence are totals of corresponding elements of the first.  We remark that the sequence $r=2s/x -1$ is called the \textit{large Schr\"{o}der} sequence \cite[p. 474]{FlSe09} (it solves $r=1+xr+xr^2$).

\ \\
\noindent {\bf (D)}  The dual interpretation of sequence expressions as set-theoretic structure specifications, is exploited in \cite{FlSaZi91,FlSe09} (influenced by \cite{Joyal81}).  One might write the set-theoretic counterpart of say, $\sfn{perm}=\sfn{set}\co\sfn{cycle}$, as $\sfn{Perm}\cong\sfn{Set}\co\sfn{Cycle}$, indicating by the initial capital letters a set-theoretic interpretation.  Essentially this is done in \cite{FlSe09}, but for brevity we only give the equations defining the counting sequences.  Table \ref{tab-seq-exps} lists some univariate examples, and table \ref{tab-biv-seqs} some bivariate ones.  These can be typed more-or-less verbatim into Haskell, as illustrated in section \ref{sec-exercising}.  Many more could be lifted from chapters I-III of \cite{FlSe09}, from the appendices of \cite{BeLaLe98}, and from \cite{Comtet,GoJa83,Stanley11,Stanley97}.
\begin{table}\[\begin{array}{rcl}
\sfn{emptySet} &=& 1\\
\sfn{singletonSet} &=& x \\
\sfn{singletonList} &=& x \\
\sfn{nonEmptyList} &=& \lst - 1 \\
\sfn{pluralList} &=& \lst - \sfn{singletonList} - 1 \\
\sfn{ordPair} &=& x^2 \\
\sfn{fibonacci} &=& \lst \co (\sfn{singletonList}+\sfn{ordPair}) \\
\sfn{cycle} &=& \log x^* \\
\sfn{oneCycle} &=& x \\
\sfn{oneOrTwoCycle} &=& \sfn{oneCycle}+x^2/2 \\
\sfn{involution} &=& \set \co \sfn{oneOrTwoCycle} \\
\sfn{nonLoopCycle} &=& \cyc-\sfn{singletonSet} \\
\sfn{derangement} &=& \set \co \sfn{nonLoopCycle} \\
\sfn{permutation} &=& \sfn{derangement} * \set \\
\sfn{nonEmptySet} &= &\set - \sfn{empty}  \\
\sfn{pluralSet} &= & \sfn{nonEmptySet} - \sfn{singletonSet} \\
\sfn{setPartition}&=& \set \co \sfn{nonEmptySet} \\
\sfn{oddNumberOfParts} &=& \sinh \co \sfn{nonEmptySet} \\
\sfn{evenSizedParts} &=& \set \co (\cosh-1) \\
\sfn{catalanTree} &=& x (\lst \co \sfn{catalanTree})\\
\sfn{cayleyTree} &=& x (\set \co \sfn{cayleyTree})\\
\sfn{connectedAcyclicGraph} &=& \sfn{cayleyTree} - \half \sfn{cayleyTree}^2 \\
\sfn{acyclicGraph} &=& \set \co \sfn{connectedAcyclicGraph} \\
\sfn{motzkinTree} &=& x(1+\sfn{motzkinTree}+\sfn{motzkinTree}^2) \\
\sfn{hipparchusSchroeder} &=& (1+x-\sqrt{1-6x+x^2})/4 \\
\sfn{largeSchroeder} &=& 2*\sfn{hipparchusSchroeder}/x -1 \\
\sfn{connectedMapping} &=& \cyc \co \sfn{cayleyTree}\\
\sfn{mapping} &=& \set \co \sfn{connectedMapping} \\
\sfn{fixedPointFree} &=& \set \co \sfn{nonLoopCycle} \co \sfn{cayleyTree} \\
\sfn{idempotent} &=& \set \co (\sfn{oneCycle} * \set) \\
\sfn{partialMapping} &=& \sfn{mapping} *(\set \co \sfn{cayleyTree})\\
\sfn{surjection} &=& \lst \co \sfn{nonEmptySet} \\
\sfn{zigzag} &=& 2(\tan + \sec) \\
bernoulli &=& x/(\exp-1)
\end{array}\]\caption{Some counting sequences \label{tab-seq-exps}}
\end{table}

Let us examine a less-than-obvious expression, \sfn{ascents} from table \ref{tab-biv-seqs}, the origin of which illustrates the principle of \textit{inclusion-exclusion} \cite{FlSe09,Zeilberger84}.  It counts permutations according to the number of ascents.  For example, the permutation $|248|3679|5|1| $
of $\{1..9\}$ has 4 up-runs demarcated with vertical bars, 3 descents, and 5 ascents.  If there are $k$ descents then there are $k+1$ up-runs.  The reversal of a permutation with $k$ descents delivers a permutation with $k$ ascents. The count $n![u^k][z^n]\sfn{ascents}$ is called an \textit{Eulerian number}, and gives the number of permutations of $\{1..n\}$ with $k$ ascents \cite{GrKnPa94}.  To come up with the sequence expression, first note that an up-run with at least one ascent corresponds to a plural set.  The counting sequence for such a set in which the $k$ of $u^k$ records the ascents is $(\sfn{pluralSet}\co uz)/u$.  Let us specify permutations in which some \textit{parts} of up-runs are identified as sets, and other elements are undistinguished:  $b(z,u)=\lst \co (z+ (\sfn{pluralSet}\co uz)/u)$.  Now propose that $\sfn{ascents}(z,u)$ is the exact counting sequence we are after, then the inclusion-exclusion principle says that $\sfn{ascents}(z,u+1)=b(z,u)$ and $\sfn{ascents}(z,u)=b(z,u-1)$, which is cited in the table.

\begin{table}\[\begin{array}{rcl}
\sfn{pascal} &=& (u+z)^* \\
\sfn{intComposition} &=& \lst \co (u*(\sfn{nonEmptyList} \co z)) \\
\sfn{schroeder} &=& z+u*(\sfn{pluralList} \co \sfn{schroeder})\\
\sfn{catalanLeaves} &=& u*z+z*(\sfn{nonEmptyList} \co \sfn{catalanLeaves}) \\
\sfn{cayleyLeaves} &=& u*z+z*(\sfn{nonEmptySet} \co \sfn{cayleyLeaves}) \\
\sfn{ebinom} &=& \set \co (z+uz) \\
\sfn{cycles} &=& \set \co (u*(\cyc \co z)) \\
\sfn{parts} &=& \set \co (u*(\sfn{nonEmptySet} \co z)) \\
\sfn{permFixedPts} &=& (\sfn{derangement}\co z)*(\set\co uz) \\
\sfn{zigzags} &=& (\sin \co u + \cos \co u)/\cos \co (u+z) \\
\sfn{ascents} &=& \lst \co (z+(\sfn{pluralSet}\co (uz-z))/(u-1))\\
\sfn{valleys} &=& \displaystyle\frac{\sqrt{1-u}}{\sqrt{1-u}-\tanh\co (z\sqrt{1-u})} \\
\sfn{powerSums} &=& \displaystyle\frac{\exp\co uz -1}{\exp\co z -1} \\
\sfn{bernoulliPoly} &=& \displaystyle\frac{z \exp \co uz}{\exp\co z -1} \\ 
\sfn{legendre} &=& (1-2uz+z^2)^{-1/2} \\
\sfn{chebyshev} &=& \displaystyle\frac{1-uz}{z^2-2uz+1} \\
\sfn{laguerre} &=& \displaystyle \frac{1}{1-z}\exp\co \frac{-uz}{1-z} \\
\sfn{hermite} &=& \exp \co (2uz-z^2)  \\
\sfn{meixner} &=& (1+z^2)^{-1/2}\exp \co (u \arctan \co z) 
\end{array}\]\caption{Some bivariate sequences \label{tab-biv-seqs}}
\end{table}

\ \\
\noindent {\bf (E)}  The sequence defined by $D\,\tan=1+\tan*\tan; \ \tan_0=0$ is the Maclaurin expansion of the tangent function (the $*$ is explicit just for emphasis).  The numbers in $t=\leo \tan=[0,1,0,2,0,16,0,272,0,7936,\ldots]$ are called \textit{tangent numbers}.  The tangent numbers count certain kinds of alternating permutations (or ordered binary trees) \cite{Stanley10}.  One can define $t$ also by $E\,t=1+t \otimes t$, where $\otimes$ is \textit{shuffle} (or \textit{Hurwitz} \cite{Keigher97}) product.
Convolution product (\ref{eq-conv-prod}) and shuffle product can be defined in head-tail form ($E \equiv\ '$):
\[\begin{array}{rcl|rcl}
(st)' &=& s't+s_0t' & (st)_0=s_0t_0 \\
(s \otimes t)' &=& s'\otimes t + s \otimes t' & (s\otimes t)_0=s_0t_0
\end{array}\]
The rule $E(s\otimes t)=E\,s\otimes t+s\otimes E\,t$ matches $D(fg)=D\,f+fD\,g$ and we have the Leibniz formulae
\[ D^n (fg) = \sum_{k=0}^n \binom{n}{k}(D^k f)  D^{n-k} g; \hspace{10mm}E^n (s \otimes t) = \sum_{k=0}^n \binom{n}{k}E^k s \otimes E^{n-k} t \]
Applying $[x^0]$ to the latter gives the pointwise definition of $\otimes$ (note that $a \otimes b=ab$ when $a$ and $b$ are scalars):
\[[x^n] (s \otimes t) = \sum_{k=0}^n \binom{n}{k}([x^0] E^k s) \otimes ([x^0]E^{n-k} t) = \sum_{k=0}^n \binom{n}{k}s_k t_{n-k} \]
A shuffle inverse, $\shi{s}{1}$, derives from the specification $s\otimes\shi{s}{1}=1$ together with the shuffle product rule (in exactly the same way that the rule for $D\,f^{-1}$ is derived):
\[ 0 = 1' = (s\otimes\shi{s}{1})' = s'\otimes \shi{s}{1} + s\otimes(\shi{s}{1})'; \ \  (\shi{s}{1})'=-s'\otimes \shi{s}{1}\otimes \shi{s}{1} \]
\noindent {\bf (F)} Let $\seqs_F(*)$ denote the ring $(\seqs_F,+,*,0,1)$, of sequences over some field $F$ (of characteristic 0) with the availability of inverses implied.  Then $(\seqs_F(*),D,\int)$ is an integro-differential algebra \cite{BuRo12}, and so too is $(\seqs_F(\otimes),E,x)$, where $x$ is the right-shift operator: $xf=x*f$.  The transform $(\leo,\loe)$ is an isomorphism between them, and is a formal Laplace transform \cite{PaEs98,Gatto12}.  In the previous example, the equation defining $\tan$ is transformed by $\leo$ into the equation defining $t$.  The sequence of factorial numbers , $x^\otimes=\leo x^*$, can be defined by applying $\leo$ to $x^*=1+xx^*$ to get $x^\otimes=1+x\otimes x^\otimes$.  From this we deduce $x^\otimes=\shi{(1-x)}{1}$.  Observe, for $n>0$: \[ [x^n]x^\otimes=[x^{n}]x\otimes x^\otimes = \sum_{k=0}^{n} \binom{n}{k}[x^k]x[x^{n-k}]x^\otimes =\binom{n}{1}[x^{n-1}]x^\otimes = n[x^{n-1}]x^\otimes \]
Also, $[x^n]x^\otimes = n[x^{n-1}]x^\otimes= (n-1)[x^{n-1}]x^\otimes+[x^{n-1}]x^\otimes=[x^n]x^2D\,x^\otimes+[x^n]xx^\otimes$ leads to the differential equation for the factorials:
\[ x^\otimes = 1 + xx^\otimes + x^2D\,x^\otimes \]
Furthermore, 
\[ (x^\otimes)'= (\shi{(1-x)}{1})'  = \shi{(1-x)}{2}; \hsp x^\otimes=1+x\shi{(1-x)}{2} \]
\noindent {\bf (G)} We recall a classic proof \cite{Niven69} of the binomial theorem.  Let $r\in F$, $\falling{r}{k}=r(r-1)\cdots (r-k+1)$ (the falling factorial), and  $\rising{r}{k}=r(r+1)\cdots (r+k-1)$ (the rising factorial): 
\[ [x^k](1+bx)^r=\frac{1}{k!}[x^0]D^k (1+bx)^r = \frac{1}{k!}[x^0]\falling{r}{k}b^k(1+bx)^{r-k} = b^k\binom{r}{k} \]
A corollary is $[x^k](x^*)^r=[x^k](1-x)^{-r}=(-1)^k\falling{(-r)}{k}/k!=\rising{r}{k}/k!=\binom{r+k-1}{r-1}$. This result, and $[x^m]\exp^n=[x^m]\exp\co nx=n^m/m!$, are basic ingredients in the search for $n$th term formulas.  They are applied next.  

\ \\
\noindent {\bf (H)} A Lagrange inversion formula \cite{Stanley99,MeSpVe06,MeSpVe07,Gessel2016} gives an expression for the $n$th term of the converse of a sequence.  For example, let $g=x(r\co g)$, then $x=g/(r\co g)=(x/r)\co g$, so $g$ is the converse of $x/r$.  Below is the Lagrange inversion formula for this case, followed by its application to the counting sequences for Catalan trees, $c=x(\sfn{list}\co c)=x(x^*\co c)$, and Cayley trees, $t=x(\sfn{set}\co t)=x(\exp\co t)$:
\begin{eqnarray*}
\,[x^n]g &=& \frac{1}{n}[x^{n-1}]r^n \\
\,[x^n]c &=& \frac{1}{n}[x^{n-1}](x^*)^n = \frac{1}{n}\binom{2n-2}{n-1} \\
\,[x^n]t &=& \frac{1}{n}[x^{n-1}]\exp^n = \frac{1}{n}\frac{n^{n-1}}{(n-1)!}=\frac{n^{n-1}}{n!}
\end{eqnarray*}
For a history of the Catalan numbers, see \cite{Pak14}.  Cayley trees are rooted versions of the connected acyclic graphs counted by Cayley in \cite{Cayley1889}.  Cayley also counted the Catalan trees in \cite{Cayley1859}, and the first part of Niven \cite{Niven69} sets out to legitimise the sequence algebra underlying Cayley's proof (Niven cites \cite{Jacobson51}, not Cayley; but Raney \cite{Raney60} cites both).

A slightly more general statement of Lagrange inversion is that it solves $h=g\co f$ (equivalently $g=h\co \conv{f}$) for $g$, where $h\in \seqs,\ f\in \seqs_C$.  The theory of Lagrange inversion sometimes employs Laurent series -- series with negative powers (or sequences with negative indicies).  In the following formula for the $n$th term of $g=h\co \conv{f}$, the coefficient of $x^{-1}$ (called the \textit{residue}) is identified:
\[ [x^n]g  =  \frac{1}{n}[x^{-1}] D\,hf^{-n}; \hsp\hsp [x^0]g  =   h_0 \]
Let $s=x(r\co s)$, $s= \conv{(x/r)}$,  and $h=g\co (x/r)$; then Lagrange inversion applied to $g=h \co s$, gives, for $h=x^k$, the $n$th term formula $[x^n]s^k = \displaystyle\frac{k}{n}[x^{n-k}]r^n$.  This result specialises, when $r=x^*$, to a variant of the \textit{cycle lemma} \cite{DeZa90}, called Raney's lemma in \cite{GrKnPa94}, which has a history in statistics \cite{Pitman98} related to \textit{Ballot numbers} \cite[p. 68]{FlSe09}.  There are various Lagrange inversion formulas and many proofs; \cite{Chen93} takes an approach that also facilitates proof of Fa\`{a} di Bruno's formula formula for $D^n(f\co g)$ \cite{Johnson02}.
 
\ \\
\noindent {\bf (I)}  The (forward) \textit{difference} operator, $\Delta s=[E-1]s=s'-s$ produces the sequence of term-to-term differences, $[x^n]\Delta s=s_{n+1}-s_n$.   The definition of an \textit{anti-difference} operator $\Sigma$ on sequences is calculated \cite{Hinze08} as a right identity to $\Delta$, with $(\Sigma s)_0=0$:
\[\Delta \Sigma s = s \Leftrightarrow  (\Sigma s)'-\Sigma s = s \Leftrightarrow  (\Sigma s)'= \Sigma s + s \Leftrightarrow  \Sigma s = 0 + x(\Sigma s + s) \Leftrightarrow \Sigma s = \frac{xs}{1-x}  \]
Thus, $\Sigma s= xx^* s$ computes all the prefix sums of $s$ (including the empty one).  Applied to $\Delta s$ we get:
\[ [x^n]\Sigma \Delta s = [x^n]xx^* \Delta s = [x^{n-1}]x^*\Delta s=\sum_{i=0}^{n-1}\left([x^{i+1}]s-[x^i]s\right) = [x^n]s-[x^0]s \]
There follows the fundamental theorem of discrete calculus (FDC) on sequences: $ s = \rep{s_0} + \Sigma \Delta s;\  s = \Delta \Sigma s$, where $\rep{a}=ax^*$ is the sequence with $a$ everywhere. 

\ \\
\noindent {\bf (J)} Here are the $E$ to $\Delta$ translations extended to powers:
\begin{eqnarray}
E^n  =  (1+\Delta)^n &=& \sum_{k=0}^{n}\binom{n}{k}\Delta^k \nonumber\\
\,[x^n]s = (E^n s)_0 & = &\sum_{k=0}^{n}\binom{n}{k}(\Delta^k s)_0  \label{eq-toDelta} \\
 \Delta^n=(E-1)^n & = &\sum_{k=0}^{n}\binom{n}{k}(-1)^{n-k}E^k \nonumber\\
 (\Delta^n s)_0 &=& \sum_{k=0}^{n}\binom{n}{k}(-1)^{n-k}s_{k} = [x^n](-x)^*\otimes s \label{eq-fromDelta}
\end{eqnarray}
The identity $\binom{n}{k}=[x^n]x^k/(1-x)^{k+1}$ turns equation (\ref{eq-toDelta}) into the Euler expansion, $s=\displaystyle \sum_k \frac{(\Delta^k s)_0 x^k}{(1-x)^{k+1}}$.  This expansion can also be derived from $s \otimes (-x)^*=\sum_k (\Delta^k s)_0x^k$, plus the facts $(-x)^*\otimes x^*=1$, and $x^k\otimes x^*=x^k/(1-x)^{k+1}$ (see \cite{BaHaPiRu17} and items {\bf K} and {\bf P}):
\begin{equation}
s=s\otimes (-x)^*\otimes x^*=(s_0+(\Delta s)_0x+(\Delta^2 s)_0x^2 + \cdots)\otimes x^*=\sum_k \frac{(\Delta^k s)_0 x^k}{(1-x)^{k+1}} \label{eq-euler-expansion}
\end{equation}
Let $g= \lne \co -x$, whence $\lne=g\co -x$, and apply Euler's expansion to $g$,
\begin{eqnarray*}
\lne &=&\left(\sum_{k}\frac{x^k}{(1-x)^{k+1}}(\Delta^k g)_0\right)\co -x=  \sum_{k}\frac{(-x)^k}{(1+x)^{k+1}}(\Delta^k g)_0\\ 
\lne(1) &=&  \sum_{k}\frac{(-1)^k}{2^{k+1}}(\Delta^k g)_0 
\end{eqnarray*} 
It is instructive to use this to approximate $\log(2)=\lne(1)$.

\ \\
\noindent {\bf (K)} The sequence ${\cal N}s=(-x)^*\otimes s=[s_0,(\Delta s)_0,(\Delta^2 s)_0, \ldots]$ is the sequence of \textit{Newton coefficients} \cite{BaHaPiRu17}.  It may also be specified by $({\cal N}s)' = {\cal N}(\Delta s);\ ({\cal N}s)_0=s_0$.  The operator ${\cal N}=((-x)^*\otimes\_)$, called the \textit{Newton transform} in \cite{BaHaPiRu17}, has the converse ${\cal N}^{-1}=(x^*\otimes\_)$, called the \textit{Binomial transform} in \cite{Hinze10}. The identity  $x^*\otimes (-x)^*=1$ holds because the head of $x^*\otimes (-x)^*$ is 1, and the tail is
\[ (x^*\otimes (1+x)^{-1})' = x^*\otimes ((1+x)^{-1}+((1+x)^{-1})') = 0 \]
The following products introduce two new rings, the \textit{Hadamard} ring $(S_R, +, -, \odot, 0, \rep{1})$, and the  \textit{infiltration} ring $(S_R,+,-,\uparrow, 0, 1)$ \cite{BaHaPiRu17}:
\[\begin{array}{lcllcl}
 (s\odot t)' &=& s' \odot t' &(s\odot t)_0 &=& s_0 t_0 \\
 (s \uparrow t)'&=&s' \uparrow t+s \uparrow t'+s' \uparrow t' &(s \uparrow t)_0&=&s_0t_0
\end{array}\]
The rules in table \ref{tab-delta-sigma-rules} apply, and $({\cal N},{\cal N}^{-1})$ is an isomorphism between the Hadamard and infiltration rings.  The following is a point-wise definition of $s \uparrow t$:
\[ [x^n]s \uparrow t=[x^n]{\cal N}({\cal N}^{-1}s \odot {\cal N}^{-1}t)=\sum_{i=0}^n\binom{n}{i}(-1)^{n-i}[x^i]((x^*\otimes s)\odot(x^*\otimes t)) \]
The proof that ${\cal N}$ is a morphism from $\odot$ to the new product $\uparrow$ can be re-imagined as a discovery of what the definition of $\uparrow$ should be.  The morphism presumption is signalled on the right below.
\[\begin{array}{lcll}
\lefteqn{({\cal N} (s\odot t))'}\\
 &=& {\cal N}(\Delta (s\odot t)) & \mbox{defn. }{\cal N}\\
&=& {\cal N}(\Delta s \odot t + s\odot \Delta t + \Delta s \odot \Delta t) & \mbox{$\Delta$-product (2)}\\
&=& {\cal N}(\Delta s \odot t) +  {\cal N}(s\odot \Delta t) +{\cal N}(\Delta s \odot \Delta t)&\mbox{morphism}\\
&=& {\cal N}\Delta s \uparrow {\cal N}t +  {\cal N}s\uparrow {\cal N}\Delta t +{\cal N}\Delta s \uparrow {\cal N}\Delta t&\mbox{morphism}\\
&=&({\cal N} s)' \uparrow {\cal N}t +  {\cal N}s\uparrow ({\cal N} t)' +({\cal N} s)' \uparrow ({\cal N} t)' &\mbox{defn. }{\cal N} \\
&=& ({\cal N} s \uparrow {\cal N} t)' & \mbox{defn. }\uparrow
\end{array}\]
\begin{table}\[\begin{array}{lrcll}
1. & \Delta (s \odot t) &=&  s \odot \Delta t + \Delta s \odot t' & \mbox{\bf $\Delta$-product (1)} \\
2. & \Delta (s \odot t) &=& \Delta s \odot t + s \odot \Delta t + \Delta s \odot \Delta t & \mbox{\bf $\Delta$-product (2)} \\
3. & \Sigma (s \odot \Delta t) &=& s\odot t - \Sigma (\Delta s \odot t') - \rep{(s\odot t)_0} & \mbox{\bf $\Sigma$-product (1)} \\

4. & \Sigma (s'\odot v) &=& s\odot (\Sigma v) - \Sigma (\Delta s \odot \Sigma v) & \mbox{\bf $\Sigma$-product (2) } \\
5. & \Sigma \Delta s \odot \Sigma \Delta u + \Sigma \Delta (s\odot u) &=& (\Sigma \Delta s)\odot u + s\odot (\Sigma \Delta u) & \mbox{\bf $\Sigma$-Baxter rule} 
\end{array}\]\caption{$\Delta-\Sigma$ rules\label{tab-delta-sigma-rules}}
\end{table}

\ \\
\noindent {\bf (L)} The following defines permutation cycle numbers, $\perms{n}{k}=n![u^kz^n]\sfn{cycles}$, and set partition numbers $\parts{n}{k}=n![u^kz^n]\sfn{parts}$.  These are also called Stirling numbers of the first and second kind, respectively.
\[\begin{array}{lcl}
\sfn{cycles} &=& \set \co (u*(\cyc \co z)) \\
\sfn{parts} &=& \set \co (u*(\sfn{nonEmptySet} \co z)) \\
\end{array}\]
The well-known recurrences \cite{GrKnPa94},
\[ \perms{n+1}{k}=\perms{n}{k-1}+n\perms{n}{k}, \hspace{10mm} \parts{n+1}{k}=\parts{n}{k-1}+k\parts{n}{k} \]
translate, using $c=\sfn{cycles}$ and $p=\sfn{parts}$, into
\[ D_z c = uc+zD_z c \hspace{20mm} D_z p = up +uD_u p \]
where $D_z$ and $D_u$ are the partial differentiation operators with respect to $z$ and $u$.  To see this, note that $\perms{n+1}{k}=(n+1)![u^kz^{n+1}]c=n![u^kz^n]D_zc$, $\perms{n}{k-1}=n![u^{k-1}z^n]c=n![u^kz^n]uc$, and so on.  The recurrences can be checked: first, $D_z c=D_z \exp\co(u \log z^*)= ucz^*= uc + uczz^*=uc+zD_z c$; and second, $up+uD_u p=up+u(\exp\co z-1)p=pu\exp\co z = D_z p$.  We may write $n![z^n]\exp\co(u\,\log\co z^*)= n![z^n](1-z)^{-u} = \rising{u}{n}$.  The cycles recurrence can also be written
$[x^k]\rising{x}{n} = [x^{k-1}]\rising{x}{n-1}+(n-1)[x^k]\rising{x}{n-1}$, which follows from $\rising{x}{n} = \rising{x}{n-1}(x+n-1)=x\rising{x}{n-1}+(n-1)\rising{x}{n-1}$.
 
\ \\
\noindent {\bf (M)} A \textit{factorial polynomial} uses falling factorials instead of powers.  For example, let $p=1+2x+x^2$, then the falling factorial counterpart is $\newton{p}=1+3\falling{x}{1}+\falling{x}{2}$.  Coefficients in $\newton{p}$ are identified by $[\falling{x}{k}]\newton{p}$, for example $[\falling{x}{1}]\newton{p}=3$.  

The symbols $\Sigma$ and $\Delta$ are overloaded as operators on factorial polynomials and obey rules identical to those for $D$ and $\int$ on polynomials: let $\newton{p}$ denote a polynomial in falling factorials, then $[\falling{x}{k}]\Delta \newton{p} = (k+1)[\falling{x}{k+1}]\newton{p}$ and $[\falling{x}{n+1}]\Sigma \newton{p} = \displaystyle\frac{1}{n+1}[\falling{x}{n}]\newton{p}$.  The fundamental theorem of the discrete calculus on factorial polynomials is immediate: $\newton{p} = \newton{p}_0 + \Sigma \Delta \newton{p};\ \newton{p}=\Delta \Sigma \newton{p}$.  In \cite{GrKnPa94}, the theorem provides one of seven ways of deducing the polynomial for summing squares: given  $\Delta \newton{p} = (1+x)^2 = 1+3\falling{x}{1}+\falling{x}{2}$, apply $\Sigma$ to both sides,  
\begin{equation}
\newton{p} - \newton{p}_0 =  1\falling{x}{1}+3/2\,\falling{x}{2} + 1/3\,\falling{x}{3};\ 
  p= 1/6\,(x+3x^2+2x^3) \label{eq-sum-squares}
\end{equation}
Note that if $p$ is a polynomial $n$th term formula for sequence $s$, $p(n)=s_n$, then $(\Delta \newton{p})(n)=(\Delta s)_n$ and $(\Sigma \Delta \newton{p})(n)=p(n)-p(0)=s_n-s_0=(\Sigma \Delta s)_n$.
 
\ \\
\noindent {\bf (N)}    
An analogue of the Maclaurin rule holds: $[\falling{x}{n}]\newton{p} = \displaystyle\frac{1}{n!}[\falling{x}{0}]\Delta^n\newton{p} = \frac{1}{n!}(\Delta^n p)(0)$.  The latter equality involves yet another interpretation of $\Delta$: $(\Delta p)(n)=p(n+1)-p(n)$.  Gregory-Newton (interpolation) formulas for $p$ of degree $m$ follow; the second (see also (\ref{eq-toDelta})) uses $\falling{n}{k}/k!=\binom{n}{k}$:
\begin{eqnarray}
\newton{p} &=& p(0)\falling{x}{0} + (\Delta p)(0)\falling{x}{1} + \frac{(\Delta^2 p)(0)}{2!}\falling{x}{2}+ \cdots+\frac{(\Delta^m p)(0)}{m!}\falling{x}{m} \label{eq-Newton-exp1} \\
p(n) &=& p(0)\binom{n}{0} + (\Delta p)(0)\binom{n}{1}+(\Delta^2 p)(0)\binom{n}{2}+\cdots +(\Delta^m p)(0)\binom{n}{m} \label{eq-Newton-exp2}
\end{eqnarray}
Let $s=[0,1,5,14,30,55, \ldots]$ be the sequence $s_n=1^2+2^2+\cdots+n^2$, for which we seek the polynomial $p$ such that $p(n)=s_n$.  Then the above expansions produce the polynomial(s) in (\ref{eq-sum-squares}).  By contrast, the Euler expansion (\ref{eq-euler-expansion}) produces the sequence expression $s=(x+x^2)/(1-x)^4$.

\ \\
\noindent {\bf (O)} We have seen $[x^k]\rising{x}{n}=\perms{n}{k}$, so we can express the polynomial for $\rising{x}{n}$ in terms of cycle numbers, and by change of signs, also the polynomial for $\falling{x}{n}$:
\[ \rising{x}{n} = \sum_{k=1}^{n}\perms{n}{k}x^k; \hspace{10mm}\falling{x}{n} = \sum_{k=1}^{n}(-1)^{n-k}\perms{n}{k}x^k; \hsp \rising{x}{0}=\falling{x}{0}=1 \]
This shows how to translate falling factorials into powers.  The converse is
\[ x^n = \sum_{k=1}^{n}\parts{n}{k}\falling{x}{k} \hspace{15mm}[\falling{x}{k}]x^n = \parts{n}{k}  \]
and for a proof see \cite[p. 343]{Brand66} and \cite[p. 262]{GrKnPa94}. 

\ \\
\noindent {\bf (P)} Infinite sequences are called  \textit{streams} when they are identified with the final object in a category of head-tail coalgebras \cite{JaRu97,Rutten05}.  This accounts for the name ``stream'' in the following:
\[\begin{array}{llllll}
\lefteqn{\mbox{\hsp Fundamental theorem of (sequence) calculus (FTC)}} \\
f & = & f_0 + \int D\, f \hspace{10mm} & f & = & D\,\int f \\
\lefteqn{\mbox{\hsp Fundamental theorem of stream calculus (FSC)}} \\
f & = & f_0 + x (E f) \hspace{10mm} & f & = & E(x f) \\
\lefteqn{\mbox{\hsp Fundamental theorem of discrete calculus (FDC)}} \\
f & = & \rep{f_0} + \Sigma \Delta\, f & f & = & \Delta\,\Sigma f\\
\end{array}\]
The co-algebraic stream calculus \cite{Rutten05} introduces a proof principle called \textit{coinduction}.  For example, the identity $x^k\otimes x^*=x^k/(1-x)^{k+1}=x^*(xx^*)^k$ used in the proof of (\ref{eq-euler-expansion}) can be proved using coinduction (it can also be proved from the point-wise definition of $\otimes$ and the binomial theorem).  Here is the gist of the coinductive proof.  Propose the relation $x^k\otimes x^*\sim x^*(xx^*)^k$.  This is used as a coinductive hypothesis.  Head-equality holds, $(x^k\otimes x^*)_0= (x^*(xx^*)^k)_0$.  The proof is completed by showing that the tails are equal under the hypothesis, signified below by the use of $(\sim)$:
\[\begin{array}{lcll}
(x^*(xx^*)^k)' &=& (x^*)'(xx^*)^k +1((xx^*)^k)'& \mbox{ conv. product rule} \\
&=& x^*(xx^*)^k + x^*(xx^*)^{k-1} & \ (x^*)'=x^* \mbox{ and } ((xf)^n)'=f(xf)^{n-1}\\
&\sim & x^k\otimes x^* +  x^{k-1}\otimes x^* & \mbox{ coinduction hyp.}\\
&=& (x^k \otimes x^*)' & \mbox{ shuffle product rule}
\end{array}\]
The bracketed (co-) in the paper's title indicates that we only touch lightly on co-algebraic concepts.  There is more to coinduction than the above example suggests, and we refer to the survey \cite{HaKuRu17} for background.  \\
\ \\
\noindent {\bf (Q)} One can check from the pointwise definitions of $D$ and $\otimes$ that $D\,f = (x\otimes f')'$.  Alternatively, equality can be proved by showing that these expressions satisfy the same head-tail equations.  The head-tail equation for $D\,f$ is calculated:
\[ D\,f=D(f_0+xf')=f'+xD\,f'=f'_0+xf''+xD\,f'=f_1+x(f''+D\,f') \]
Hence, $(D\,f)_0=f_1$ and $(D\,f)'=f''+D\,f'$.   Now let $F\,f=(x\otimes f')'$.  We find $(F\,f)_0=f_1$, and 
$(F\,f)' = (f'+(x\otimes f''))' = f'' + F\,f'$.  Thus, $D\,f=F\,f$ since they satisfy the same head-tail equations.

To give a little more feeling for the coinduction game, let us reveal the machine-level minutiae that proves $D(f\co g)=(D\,f\co g)D\,g$.  We use head-tail properties such as $(f\co g)_0=f_0;\ (f\co g)' = (f'\co g)g'$ (see section \ref{sec-Haskell}, equation (\ref{eq-comp})).  We will also make use of $(D h)_0=h'_0$.
Head equality is confirmed:
\[(D(f\co g))_0  =   (f\co g)'_0 = (f' \co g)_0 g'_0 =f'_0 g'_0 = (D f \co g)_0 (D g)_0 = ((D f  \co g)D g)_0\]
Tails are proved equal under the coinductive hypothesis, $D(f \co g) \sim  (D f \co g)D g$:
\[\hspace{-10mm}\begin{array}{lcll}
\lefteqn{((D f \co g)D g)'} \\
& = & (D f\co g)'D g + (Df\co g)_0(Dg)' &(')\mbox{-product}\\
 & = & ((Df)'\co g)g'Dg+(f'\co g)_0(Dg)' &(')\mbox{-composition}\\
 & = & (f''\co g)g'Dg+(Df'\co g)g'Dg+(f'\co g)_0(Dg)' &D\mbox{-defn and $\co$-distr}\\
 & \sim & (f'\co g)'Dg + D(f'\co g)g' +(f'\co g)_0(Dg)' & (')\mbox{-comp., coinduction}\\
 & = & (f'\co g)'(g'+xDg')+ D(f'\co g)g' + (f'\co g)_0(g''+Dg') & \mbox{expand }D\,g \mbox{ and } (D\,g)'\\
  & =& \lefteqn{D(f'\co g)g' + (f'\co g)'xDg'+(f'\co g)_0Dg' +(f'\co g)'g'+(f'\co g)_0g''}\\
  & = & D(f'\co g)g' + (f'\co g)Dg' + ((f'\co g)g')' & \mbox{head-tail, }h=h_0+xh'\\
  & = & D((f'\co g)g') +((f'\co g)g')' & D\mbox{-product}\\
  & = & D(f\co g)' + (f\co g)'' & (')\mbox{-composition, twice}\\
  & = & (D(f\co g))' &D\mbox{-defn}
\end{array}\]
 
\ \\
\noindent {\bf (R)} Coinduction is also used in \cite{Rutten03a} to show how continued fractions can be obtained from combinatorially-inspired automata.  For example, the tangent sequence can be defined by $t=xu_1$, where $u_k = 1/(1-k(k+1)x^2u_{k+1})$.  Thus $t$ can be displayed as a continued fraction:
\[ t = \frac{x}{1-\frac{\displaystyle 1* 2x^2}{\displaystyle 1-\frac{2* 3x^2}{\displaystyle 1-\frac{ 3* 4x^2}{\ddots}}}} \]
More combinatorially-inspired continued fraction expressions for sequences appear in \cite{Flajolet80,GoJa83,Rutten03a}.  Below is one for $\leo_z\sfn{cycles}=\shi{(1-z)}{u}$ ($\leo_z$ removes the factorial divisors of powers of $z$).
\[ \shi{(1-z)}{u} = \frac{1}{1-uz-\frac{\displaystyle 1uz^2}{\displaystyle 1-(u+2)z-\frac{2(u+1)z^2}{\displaystyle 1-(u+4)z-\frac{ 3(u+2)z^2}{\ddots}}}} \]
Setting $u=1$ and $z=x$ gives a continued fraction for the factorials, $\shi{(1-x)}{1}=x^\otimes$.
 
\ \\
\noindent {\bf (S)} A $k$th-order linear ordinary homogeneous differential equation,
\[ b_kD^k f + b_{k-1} D^{k-1}f + \cdots + b_0 f =0 \] 
can be written $\lde{b}(D)f = 0$.  Similarly, a difference equation (also called a recurrence equation) can be written $\lde{b}(E)s = 0$.   Let $\revp{b}=b_k+b_{k-1}x^1+\cdots+b_1x^{k-1}+b_0x^k$ be the \textit{reverse} of $b=b_0+b_1x+b_2x^2+\cdots + b_kx^k$.  Klarner \cite{Klarner69} presents this fact: the solution $s$ to $\lde{b}(E)s = 0$ is
\[ s= \frac{(\revp{b}*\mbox{inits})[0..k-1]}{\revp{b}} \]
where inits$=s_0+ s_1x+\cdots + s_{k-1}x^{k-1}=s[0..k-1]$ are the initial $k$ elements of $s$.  Also
\[ \lde{b}(E)s=0 \Leftrightarrow \lde{b}(D)\loe s = 0 \]
For example, $E^2s+s=0,\ s_0=0,\ s_1=1$ has solution $x/(1+x^2)$, and $\loe s=\sin$, the Maclaurin expansion for $\sin$, that is, the solution to $D^2s+s=0;\ s_0=0,\ s_1=1$.  Another example is $[z^n]C-2u[z^{n-1}]C+[z^{n-2}]C=0,\ C_0=1,\ C_1=u$, where $[z^n]C$ is a Chebyshev polynomial \cite{Brand66} in $u$.  Then, $b=1-2uz+z^2=\revp{b}$, and 
\[ C(z,u) = \frac{(\revp{b}*(1+uz))[0..1]}{\revp{b}}=\frac{1-uz}{1-2uz+z^2} \]
By the translation rules of item {\bf J}, difference equations can be written using either $E$ or $\Delta$.   
The equation $\lde{b}(E)s = 0$ transforms into $\lde{b}(1+\Delta)s = 0$, or $\hat{b}(\Delta)s = 0$, where $\hat{b}=b\co (1+x)$.  The converse is $b=\hat{b}\co (x-1)$, reflecting $\hat{b}(\Delta)=\hat{b}(E-1)$.  Clearly, $b=b\co(1+x)\co(x-1)$.

\ \\
\noindent {\bf (T)} A sequence $s$ is called \textit{rational} if it is the quotient, $s=a/b$, of polynomials; it is called a \textit{LODE solution}, written LODE($s$), if it is a solution of a linear ordinary homogeneous difference equation; and it is called \textit{recognizable} if it is the behaviour of a finite automaton.  Then,
\[ \mbox{rational}(s) \Leftrightarrow \mbox{LODE}(s) \Leftrightarrow \mbox{recognizable}(s) \]
Following \cite{Eil:74}, a finite automaton can be modelled as a system of linear equations over sequences, $E\,S=A S;\ S(0)=v$.  Here, matrix $A$ records transition labels connecting pairs of states, and $S$ is a vector of sequences, one for each state, with initial values $S_i(0)=v_i$.  The solution is $S=(Ax)^*v$ where
\[ (Ax)^* = I+Ax+A^2x^2+A^3x^3+\cdots = \sum_{i\geq 0}A^ix^i \]
is the Kleene star.  Notice that $(Ax)^*$ can be viewed as a matrix of sequences or as a sequence of matrices.  Using $(Ax)^*=(I-Ax)^{-1}$, we get $(I-Ax)S=v$ and Cramer's rule applies: $S_i=\displaystyle\frac{\det (I-Ax)[i\leftarrow v]}{\det(I-Ax)}$ (column $i$ replaced by $v$).  Thus $S_i$ is rational.  Justification of the above equivalences is completed by noting that a LODE can be transformed into a system of linear equations.  We remark that a quotient of polynomials, $a/b$, can also be written as the solution to a system of linear equations \cite{Rutten05}.

\ \\
\noindent {\bf (U)}  Let $\lde{b} = \det (xI-A) = b_0 + b_1x^{1} + \cdots +b_{n-1}x^{n-1}+ x^n$ be the characteristic polynomial of matrix $A$.  The Cayley-Hamilton theorem \cite{Eil:74,Klarner76} can be stated as $\lde{b}(A)=0$, or as  $\lde{b}(E)(Ax)^* = 0$.  This will hold if, taking the sequence of matrix powers, $(Ax)^*$, now as a matrix of sequences, we have $\lde{b}(E)((Ax)^*)_{ij} = 0$, which in turn holds if
$((Ax)^*)_{ij} = \displaystyle\frac{a}{\revp{b}}$.  We have $\revp{b}=x^n\lde{b}(1/x)=\det (I-xA)$.  So we are done if we come up with an $a$ such that
\[ ((Ax)^*)_{ij} = \frac{a}{\det (I-Ax)} \]
Let $M=Ax$, and $J=M^*\col{\uvec{j}}= \col{\uvec{j}}+MJ$, where $\uvec{j}$ is the vector with 1 at position $j$ and zero elsewhere.  Then, $(I-M)J=\col{\uvec{j}}$ and Cramer's rule delivers the $a$ we are looking for:
\[ (M^*)_{ij}=J_i = \frac{\det (I-M)[i\leftarrow\uvec{j}]}{\det (I-M)} \]
 
\ \\
\noindent {\bf (V)}  Consider the matrix exponential, $\exp \co Ax=\loe (Ax)^*$ as a matrix of sequences.  We know that $(Ax)^*$ solves $E\,S=AS,\ S(0)=I$, and $(Ax)^*[0..k-1]=[I,A,A^2,\ldots,A^{k-1}]$.  Cayley-Hamilton says that $\lde{b}(E)(Ax)^* =0$, where $\lde{b}$ is the characteristic polynomial of $A$, of degree $k$, say.  By uniqueness of solution, we have \cite{Leonard96,Liz98, Kwapisz98} $\phi=(Ax)^*$ if $\lde{b}(E)\phi =0$ and $\phi[0..k-1]=[I,A,A^2,\ldots,A^{k-1}]$.  Let $S$ be a vector of sequences such that $\lde{b}(E)S_i=0$ and  $S_i[0..k-1]=x^i$.  Set
\[ \phi=S_0I+S_1A+\cdots S_{k-1}A^{k-1} = \sum_{i=0}^{k-1}S_iA^i\]
Clearly $\phi[0..k-1]=[I,A,A^2,\ldots,A^{k-1}]$, and $\lde{b}(E)\phi = 0$ because
\[  \sum_{j=0}^{k}\lde{b}_jE^j\phi = \sum_{j=0}^{k}\lde{b}_j\left(\sum_{i=0}^{k-1}E^jS_iA^i\right) =\sum_{i=0}^{k-1}\left(\sum_{j=0}^{k}\lde{b}_jE^jS_i\right)A^i =\sum_{i=0}^{k-1}(\lde{b}(E)S_i)A^i = 0 \]

\ \\
\noindent {\bf (W)} The elements of the sequence $B=x/(\exp-1)=[1,-1/2,1/6,0,-1/30,0,1/42,0,\ldots]$ are called \textit{Bernoulli numbers}. The corresponding recurrence is calculated from $B\,\exp = B+x$:
\[ n![x^n]B\,\exp = n![x^n](B+x); \hspace{10mm} \sum_{k=0}^{n} \binom{n}{k} B_k = B_n + [n=1] \]
Bernoulli numbers are used by Graham \textit{et al} \cite{GrKnPa94} for the most impressive of their deductions of the polynomial that sums squares -- impressive because it defines the formulas for all powers at once. Let $S_{(n)}$ be the sequence such that $m![x^m]S_{(n)}$ is the sum of the $m$th powers of the naturals to $n-1$.  Then
\begin{eqnarray*}
 m![x^m]S_{(n)} & = & \sum_{k=0}^{n-1} k^m\ =\ \sum_{k=0}^{n-1} m![x^m]\exp \co kx  \\
 S_{(n)} & = & \sum_{k=0}^{n-1} \exp^k  \ =\ \frac{\exp^n - 1}{\exp -1}\ =\ \frac{\exp \co nx - 1}{\exp -1} =  B\,\frac{\exp \co nx -1}{x} \\
m![x^m]S_{(n)} & = & m!\sum_{k=0}^m \frac{B_{m-k}}{(m-k)!} \frac{n^{k+1}}{(k+1)!} = \sum_{k=0}^m \binom{m}{k}B_{m-k} \frac{n^{k+1}}{k+1}
\end{eqnarray*}
Now replace $n$ by $u$ and $x$ by $z$ to get the expression advertised in the introduction:
\[ S = \frac{\exp \co uz - 1}{\exp\co z -1} \]
Then $m![z^m]S$ is a polynomial of degree $m+1$ in $u$, and $(m![z^m]S)(n)=m![x^m]S_{(n)}$.

\ \\
\noindent {\bf (X)} Observe that $B_1=-1/2$ is non-zero whilst all the other odd-degree coefficients of $B$ appear to be zero.  Perhaps if we make $B_1$ zero then we will have a sequence which can be \textit{proved} to be even (i.e. with zeros at odd positions).  Adding $\displaystyle\half x$ to $B$ cancels $B_1$:
\[ C = B+\frac{x}{2}=\frac{x}{\exp - 1}+\frac{x}{2}=\frac{x}{2}\,\frac{\exp+1}{\exp-1} \]
Recall $\displaystyle \coth = \frac{\exp+\exp\co-x}{\exp-\exp\co-x}$, so $C = \displaystyle\frac{x}{2} (\coth \co \frac{x}{2})$, from which eveness, $C\co -x = C$, can be deduced.  Thus,
\[ C = B+\frac{x}{2} = \sum_k \frac{B_{2k}}{(2k)!}x^{2k}\]
Now $C\co 2x = x\coth$, and, using $x\cot = ix(\coth\co ix)$ and $\tan = \cot - 2\cot\co 2x$, we get
\begin{eqnarray}
x\cot = C\co 2ix &=&\sum_k (-1)^k2^{2k}\frac{B_{2k}}{(2k)!}x^{2k}\label{eq-xcot1} \\
\tan &=& \sum_k (-1)^{k-1}4^k(4^k-1)\frac{B_{2k}}{(2k)!}x^{2k-1} 
\end{eqnarray}
With a bit of analysis (reals, $\cot(x)$ an analytic function with period $\pi$, and uniqueness of series expansion), one can deduce another series for $x\cot$, due to Euler.  The omitted analysis \cite{Knopp51,AiZi04} is hidden in the first equals sign:  
\[ x\cot(x) = 1-2\sum_{n=1}^\infty \frac{x^2}{n^2\pi^2-x^2} = 1-2\sum_{n=1}^\infty \frac{x^2}{n^2\pi^2}\,\left(\frac{x^2}{n^2\pi^2}\right)^* = 1-2\sum_{k=1}^\infty \frac{x^{2k}}{\pi^{2k}}\sum_{n=1}^\infty \frac{1}{n^{2k}}\]
Equating coefficients with those in the expansion (\ref{eq-xcot1}) yields, for $k>0$,
\[ [x^{2k}]x\cot = \frac{-2}{\pi^{2k}}\sum_{n=1}^\infty \frac{1}{n^{2k}} = (-1)^k2^{2k}\frac{B_{2k}}{(2k)!} \]
Therefore, the values of Riemann's $\zeta(s)=\displaystyle\sum_{n\geq 1}\frac{1}{n^s}$ at even positive integers is given by
\[ \zeta(2k)= (-1)^{k-1}2^{2k-1}\frac{B_{2k}}{(2k)!}\pi^{2k} \]

\noindent {\bf (Y)} 
The Formal Taylor Theorem may be expressed:
\[ f\co(u+z) = f\co u+((D\,f)\co u)z+\frac{(D^2 f)\co u}{2!}z^2+ \frac{(D^3 f)\co u}{3!}z^3+ \cdots \]
Write $f\co (u+z) = g_0+g_1z+g_2z^2+g_3z^3 + \cdots$.  Let $z=0$, then $g_0=f(u)$.  Differentiate with respect to $z$: $D_z (f\co(u+z)) = g_1+2g_2z+3g_3z^2+ \cdots$.  Note that $D_z (f\co(u+z)) = f_1+2f_2(u+z)+3f_3(u+z)^2+\cdots=(D\,f)\co (u+z)$.  Let $z=0$, then $g_1=(D\,f)\co u$.  Differentiate again: $D^2(f\co (u+z)) = 2g_2+3!g_3z + \cdots$.  Let $z=0$, then $g_2=((D^2 f)\co u)/2$. And so on. The Maclaurin expansion is the special case with $u=0$, and the Taylor expansion of $x^n\co(u+z)$ is an instance of the binomial theorem.  Lipson \cite{Lipson81} uses the theorem in the application of Newton's iterative root-finding algorithm to polynomial equations over sequences (see also \cite{PiSaSo12}).  

\ \\
\noindent {\bf (Z)} 
The following manipulations, originating with Lagrange, have a captivating charm (even if they lack rigour).  We adapt them from \cite{GrKnPa94} to show that elementary sequence algebra plays a role through to the final chapter of that book (where, however, things become more demanding).  In the Formal Taylor theorem, let $f$ be a polynomial, $z=1$, $u=x$, and employ an operator style:
\begin{eqnarray*}
E f &=& f\co (x+1)= f+(D\,f)+\frac{D^2 f}{2!}+ \frac{D^3 f}{3!}+ \cdots =[1+D+\frac{D^2}{2!}+ \frac{D^3}{3!}+ \cdots]f \\
& =& \exp(D)f
\end{eqnarray*}
Putting $\Delta = E-1=\exp(D) -1$ together with $\Delta \Sigma f= f$, suggests $\Sigma =(\exp(D) -1)^{-1}$.  Then $B=x/(\exp-1)$ applied to $D$ is $D\Sigma$, so $\Sigma = D^{-1}B(D)$.  Expanding this, and writing $\int$ for the first term $D^{-1}$ (since $B_0=1$), gives a ``template'' version of the Euler-Maclaurin summation formula \cite{GrKnPa94,Brand66}.
\[ \sum = \int + \sum_{k\geq 1} \frac{B_k}{k!}D^{k-1} \]
Now introduce limits: 
\begin{equation}
\textstyle\sum_a^b\, f = \displaystyle\int_a^b f + \left.\sum_{k\geq 1} \frac{B_k}{k!}D^{k-1}f\,\right|_a^b \label{eq-Euler-Maclaurin}
\end{equation}
An application to $x^2$ gives yet another derivation \cite{GrKnPa94} of the sum-of-squares formula:
\begin{eqnarray*}
\textstyle\sum_0^n\, x^2 &=& \frac{1}{3}n^3 + \left.\left(-\frac{1}{2}x^2+\frac{1}{12}2x\right)\right|_0^n \\
& =& \frac{1}{3}n^3 -\frac{1}{2}n^2 +\frac{1}{6}n
\end{eqnarray*}
The way the limits appear on the summation sign has significance: $\textstyle\sum_0^n\, x^2=\displaystyle\sum_{x=0}^{n-1}x^2$.  The definite summation symbol follows the pattern of definite integration:
\[\begin{array}{rcccccccl}
g = D f &\Rightarrow&\displaystyle\int_a^b g &=& \left.\mbox{\rule{0mm}{2.5ex}}f\,\right|_a^b &=& f(b)-f(a) \\
\mbox{\rule{0mm}{4ex}}g= \Delta f&\Rightarrow&\textstyle\sum_a^b\, g &=& \left.\mbox{\rule{0mm}{2.5ex}}f\,\right|_a^b &=& f(b)-f(a) &=&\displaystyle\sum_{x=a}^{b-1}f(x+1)-f(x)\\
 &&&& &=& \displaystyle\sum_{x=a}^{b-1}\Delta f(x) &=&\displaystyle\sum_{x=a}^{b-1} g(x)
\end{array}\]
Let's add $f(b)$ to both sides of (\ref{eq-Euler-Maclaurin}) and separate out $B_1=-1/2$:
\begin{equation}
\sum_{x=a}^{b}\, f = \displaystyle\int_a^b f + \frac{1}{2}(f(b)+f(a)) + \left.\sum_{k\geq 2} \frac{B_k}{k!}D^{k-1}f\,\right|_a^b \label{eq-Euler-Maclaurin2}
\end{equation}  
The Euler-Maclaurin formula can also be applied to non-polynomial functions.  Let us illustrate this, without justification.  To compute $\displaystyle\sum_{x=1}^\infty 1/x^2$, set $S_9=\displaystyle\sum_{x=1}^9 \frac{1}{x^2}=1.5397677310$ and then apply (\ref{eq-Euler-Maclaurin2}) to $g=1/(x+10)^2$:
\[\sum_{x=1}^\infty \frac{1}{x^2} = S_9+\sum_{x=0}^\infty g(x)= S_9+\int_0^\infty g +\frac{1}{2}g(0) + \left.\sum_{k\geq 2} \frac{B_k}{k!}D^{k-1}g\,\right|_0^\infty \]
Applying the formula up to $B_4$ gives $\zeta(2)=\pi^2/6\approx 1.64493407$.
\section{A programming delight}\label{sec-Haskell}
McIlroy \cite{McIlroy99,McIlroy01}, influenced by \cite{Kar97} and others, has gifted us some ``tiny gems'' of program definitions for implementing sequence manipulations.  The definitions are written in Haskell, and are effortlessly derived by mathematical reasoning.  A textbook introduction appears in \cite{DoVE12}.  We want to entice the reader to type up and experiment with the Haskell code (but a down-loadable file is available).  The code has been tested in the Haskell GHCi system, and also in the Hugs98 system (an older system, but well-suited to beginners).  Both  GHCi and Hugs98 are freely available on the web at www.Haskell.org.  

In this section we present all of the definitions, thus duplicating some of the contents of \cite{McIlroy99,McIlroy01}; however, there are modifications and additions.  The fact that the definitions have no pre-requisites, other than the standard Haskell Prelude, means that one can take ``deep'' ownership, building things from the ground up.  This contrasts to using a sophisticated computer algebra system -- something perhaps for the newcomer to move on to with greater appreciation.

Haskell \cite{Haskell} has evolved to be a fairly large and sophisticated language, but we shall stick to a modest subset.  It is expected that the reader can comprehend Haskell from examples.  The language gives types to objects and variables, and within context, the most general type is used.  Haskell's lists are used to represent sequences (some may prefer to introduce a new type for sequences, but that introduces an overhead which we want to avoid).  In Haskell, head-tail decomposition $s_0+xs'$ becomes {\tt s0:s'}.  Here are some list-processing functions:
\begin{verbatim}
    take n _ | n<=0     = []
    take _ []           = []
    take n (s0:s')      = s0: take (n-1) s'
    map f []            = []
    map f (s0:s')       = f s0 : map f s'
    iterate f z         = z: iterate f (f z)
    foldr f z []        = z
    foldr f z (s0:s')   = f s0 (foldr f z s') 
    scanl op q s        = q: (case s of 
                                []    -> []
                                s0:s' -> scanl op (op q s0) s')
    zip (s0:s') (t0:t') = (s0, t0): zip s' t'
    zip _ _             = []
    zipWith op s t      = [op sn tn | (sn,tn) <- zip s t] 
\end{verbatim} 
These definitions implement the following functions which feature in the algebra of program calculation \cite{Bir:89,BiDM97}:
\[\begin{array}{lcl}
\sfn{take}\ n\ s &=& [s_0, s_1, \ldots, s_{n-1}] \\
\sfn{map}\ f\ s &=& [f s_0, f s_1, f s_2, \ldots ] \\
\sfn{iterate}\ f\ z &=& [z, f\ z, f (f z), f^3 z, \ldots ] \\
\sfn{foldr}\ f\ z\ s &=& f\ s0\ (f\ s1\ (f\ s2\ ( \ldots\ z)\ldots ) \\
\sfn{scanl} \oplus q\ s &=& [q, q \oplus s_0, (q \oplus s_0) \oplus s_1, ((q \oplus s_0) \oplus s_1) \oplus s_2, \ldots ] \\
\sfn{zip}\ s\ t &=& [(s_0,t_0),(s_1,t_1),(s_2,t_2),\ldots ] \\
\sfn{zipWith} \odot s\ t &=& [s_0 \odot t_0, s_1 \odot t_1, s_2 \odot t_2, \ldots] 
\end{array}\]
The types deduced are
\begin{verbatim}
    take    :: Int -> [a] -> [a]
    map     :: (a -> b) -> [a] -> [b]
    iterate :: (a -> a) -> a -> [a]
    foldr   :: (a -> b -> b) -> b -> [a] -> b 
    scanl   :: (a -> b -> a) -> a -> [b] -> [a]
    zip     :: [a] -> [b] -> [(a,b)]
    zipWith :: (a -> b -> c) -> [a] -> [b] -> [c]
\end{verbatim}
Type expressions are built from type names such as {\tt Int}, type variables such as {\tt a}, and type constructors such as {\tt ->} (which associates to the right).  Two more examples of type constructors are: {\tt [a]} is the type for lists of objects of type {\tt a}, and {\tt [(a,b)]} is the type for lists of pairs of objects.  Clearly, functions can take functions as arguments, in which case they are called \textit{higher-order}.  Lazy evaluation is used, so that for example, {\tt scanl} will produce the first element, {\tt q}, of the result without needing to know anything about its list argument, {\tt s}.  The definition of {\tt zipWith} illustrates the so-called \textit{list comprehension}.  An alternative definition uses {\tt map} and {\tt uncurry}:
\begin{verbatim}
    uncurry        :: (a -> b -> c) -> (a,b) -> c
    uncurry op p   = op (fst p) (snd p)
    zipWith op s t = map (uncurry op) (zip s t)
\end{verbatim}
The partner to {\tt uncurry} is {\tt curry f x y = f (x,y)}.  An alternative definition illustrates use of a lambda expression: {\tt curry f = \verb+\+x y -> f (x,y)}.  The reader may like to supply the type.  These two functions are so-named because Haskell Curry was an early advocate of the associated equivalence \cite{HuHuPJWa07}.  The Haskell Standard Prelude defines {\tt zipWith} without using {\tt zip}, and then defines {\tt zip = zipWith (,)}.

The following specifies that a type {\tt a} is classified as {\tt Num} if it has the operations (or methods) listed here:  
\begin{verbatim}
class Num a where
    (+), (-), (*)        :: a -> a -> a
    negate, abs, signum  :: a -> a
    fromInteger          :: Integer -> a
    x - y                = x + negate y
    negate x             = 0 - x
\end{verbatim}
All of the foregoing definitions are in the Standard Prelude.  From here on, the code needs to be supplied.  The program file starts with a few specified lines: the first line hides the Prelude definition of {\tt cycle} because we are going to re-define it for other purposes; the second line says we need rational numbers; the third line gives the order in which to resolve ambiguity in numerical data.
\begin{verbatim}
    import Prelude hiding (cycle)
    import Data.Ratio
    default (Integer, Rational, Double)
\end{verbatim}  
We start by declaring how sequences, {\tt [a]}, become an instance of {\tt Num}.  A prerequisite is that {\tt a} is an instance of {\tt Eq} and {\tt Num}, indicated by {\tt (Eq a, Num a) =>}.  The definition of {\tt (-)} is derived from {\tt negate}.  
\begin{verbatim}
    instance (Eq a, Num a) => Num [a] where
      negate            = map negate 
      f+[]              = f
      []+g              = g
      (f0:f')+(g0:g')   = f0+g0 : f' + g'
      []*_              = []
      (0:f')*g          = 0 : f'*g
      _*[]              = []
      (f0:f')*g@(g0:g') = f0*g0 : (f0*|g' + f'*g)
      fromInteger c     = [fromInteger c]
      abs _             = error "abs not defined on sequences"
      signum _          = error "signum not defined on sequences"
\end{verbatim}
Observe that addition is not defined by {\tt f+g = zipWith (+) f g} (why?).  Convolution product is derived from (\ref{eq-conv-prod}), but there are some things to note.  Firstly, if $f_0=0$ then the zeroth term of the result is 0 and is delivered immediately.  This may be regarded as a controversial quirk, but it enables certain equations to be used directly, as in the following (Catalan) example -- in examples, definitions (which are placed in a program file) are interspersed with interactive requests for expression evaluation, indicated by the prompt ``\verb+> +''. 
\begin{verbatim}
    x :: Num a => [a]
    x = [0,1]
    b = 1 + x*b^2
    > take 8 b
    [1,1,2,5,14,42,132,429]
\end{verbatim}
Secondly, the notation {\tt g@}, is read as ``{\tt g} as''.  Thirdly, there are clauses for finite sequences -- the empty list behaves here like zero (but note that 0 is embedded as {\tt [0]}).  Fourthly, the term $f_0g'=[f_0]g'$ becomes an explicit scalar product using \verb+*|+, which is defined as an infix operator with precedence 7 (the same as {\tt *}, and higher than {\tt +}).  The definition contains {\tt (a*)}, illustrating  the creation of a function by partial application of an operator (called \textit{sectioning}).
\begin{verbatim}
    infix 7 *|
    (*|) :: Num a => a -> [a] -> [a]
    a *| f = map (a*) f
\end{verbatim}
A function like {\tt map} is said to be \textit{polymorphic} because any type can be assigned to its type variables (subject to consistency).  By contrast, scalar multiplication {\tt (*|)} has a \textit{qualified} (\textit{constrained} or \textit{parametric}) polymorphic type: its type variable can range over only instances of class {\tt Num}.  The type stated could be omitted because it can be inferred due to the presence of {\tt *}.  On the other hand, if the explicit type given to {\tt x} above was omitted, then Haskell would infer {\tt x::[Integer]} and this is a \textit{monomorphic} type which would restrict the use of {\tt x}.  With the qualified polymorphic type, {\tt x} can appear in an expression where a sequence of elements of type {\tt N} is expected, as long as {\tt N} is an instance of {\tt Num}.  Any instance {\tt N} of {\tt Num} must provide a {\tt fromInteger} method that shows how to embed integers into {\tt N}, so {\tt x} would be interpreted as {\tt [N.fromInteger 0, N.fromInteger 1]}.

The {\tt Num} class invites comparison with the specification of the signature of a ring.  Likewise, Haskell's {\tt Fractional} class may be compared to a ring-with-division, because a (partial) division operator {\tt (/)}, or a multiplicative inverse ({\tt recip}), is required.  Rational numbers form the archetypal instance of {\tt Fractional}, and any instance {\tt F} must show how to embed the rationals in {\tt F} by defining {\tt fromRational :: Rational -> F}.  Division on sequences, $f/g$, requires calculating the quotient $q$ satisfying $f=qg$: 
\begin{eqnarray*}
  f_0 + xf' & = & (q_0 + xq')g = q_0g + xq'g = q_0g_0 + x(q_0g'+q'g) \\
  	f_0  =  q_0g_0 &;& f'  =  q_0g'+q'g \\
  	q_0  =  f_0/g_0 &;& q'  = (f'-q_0g')/g \\
  	q & = & f_0/g_0 + x(f'-q_0g')/g
\end{eqnarray*}
Now we can say, at least approximately, how sequences become a ring-with-division:
\begin{verbatim}
    instance (Eq a, Fractional a) => Fractional [a] where
      recip f           = 1/f
      _/[]              = error "divide by zero."
      []/_              = []
      (0:f')/(0:g')     = f'/g'
      (_:f')/(0:g')     = error "divide by zero"
      (f0:f')/g@(g0:g') = let q0=f0/g0 in q0:((f' - q0*|g')/g)
      fromRational c    = [fromRational c]
\end{verbatim}
These simple definitions confront us with some of the difficulties in coding a satisfactory division operation that works for both finite and infinite sequences.  One should investigate questions like: are $f/g=f*(1/g)$ and $f/f=1$ faithfully implemented?  To keep things simple, compromises have to be made.

Arithmetic and the convolution product rule are used to calculate a head-tail definition for square root, $\sqrt{f}$.  The starting point is $f=\sqrt{f}\sqrt{f}$, and we calculate $(\sqrt{f})_0$ and $\sqrt{f}'$:
\begin{eqnarray*}
 f_0 = (\sqrt{f}\sqrt{f})_0 &=& (\sqrt{f})_0(\sqrt{f})_0\\
(\sqrt{f})_0 &=& \sqrt{f_0} \\
f'=(\sqrt{f}\sqrt{f})' & = & \sqrt{f_0}\sqrt{f}'+\sqrt{f}'\sqrt{f} \\
\sqrt{f}' & = & f'/(\sqrt{f_0}+\sqrt{f}) \hsp (f_0\neq 0)
\end{eqnarray*}
We shall trivialise $\sqrt{f_0}$ and restrict square root to fractional sequences with constant term 1.   In the following code, the first clause is suggested by the identity $\sqrt{x^2f}=x\sqrt{f}$, and the more general $\sqrt{x^{2n}f}=x^n\sqrt{f}$ is handled by recursion.  An alternative definition of square root is derived in \cite{McIlroy01} by differentiating $r^2=f$, rearranging and then integrating (which the reader may like to try).
\begin{verbatim}
    sqroot (0:0:f'') = 0:sqroot f''
    sqroot f@(1:f')  = 1:(f'/(1+sqroot f))
\end{verbatim}
Sequence composition, $f \co g  = \displaystyle \sum_{n} f_n g^n$, is expanded thus:
\begin{eqnarray*}
f \co g & = & f_0 + f_1g^1 + f_2g^2 + f_3g^3 + \cdots \\
		&		= & f_0 + g (\tail{f} \co g) =  f_0 + (g_0 + x\tail{g}) (\tail{f} \co g) \\
		&		= & f_0 + g_0(\tail{f} \co g) + x\tail{g}(\tail{f}\co g) 
\end{eqnarray*}
When $f$ is infinite, $g_0(\tail{f} \co g)$ is not computable unless $g_0=0$.  When $g_0=0$ we get
\begin{equation} (f\co g)_0=f_0; \hspace{10mm} (f\co g)'=g'(f'\co g) \label{eq-comp} \end{equation}  
However, $f\co g$ is computable for $g_0\neq 0$ when $f$ is finite ($p\co [a]=[p(a)]$, $p$ a polynomial).  So we admit a potentially non-terminating clause and it is up to us to use it with care:
\begin{verbatim}
    [] `o` _              = []
    (f0:f') `o` g@(0:g')  = f0: g'*(f' `o` g)
    (f0:f') `o` g@(g0:g') = [f0] + (g0*|(f' `o` g))+ 
                             (0:g'*(f' `o` g))
\end{verbatim}
The definition $f \co g  =  \sum_{n} f_n g^n $ reveals $x$ to be a left and right identity of composition.  Composition distributes leftwards through sum, product, and quotient.

\noindent To calculate the converse, $g=\conv{f}$,  expand the composition $f \co g = x$:
\[ f_0 + x\tail{g}(\tail{f} \co g) = x = 0+x1\]
Hence, $\tail{g}(\tail{f} \co g)=1$, and 
\[ g_0=0;\hsp\hsp g'=1/(f'\co g) \]
The program code is: 
\begin{verbatim}
    converse(0:f') = g where g = 0: 1/(f' `o` g)
\end{verbatim}
For the reciprocal of $\tail{f} \co g$ to be defined, it is necessary that $\tail{f}_0$ is invertible, which entails that $\tail{f}$ is invertible.  The set of such ``conversible'' sequences forms a group, $(\seqs_C,\co, \conv{()},x)$. 

\noindent The transforms $\leo$ and $\loe$ are given names {\tt e2o} and {\tt o2e}, respectively \cite{DoVE12}.  Here they are, together with some other useful sequences:
\begin{verbatim}
    e2o f    = zipWith (*) f facs 
    o2e f    = zipWith (/) f facs
    from     :: Num a=>a->[a]
    from     = iterate (+1)
    nats, pos, zeros, facs :: Num a=>[a]
    nats     = from 0
    pos      = from 1
    zeros    = 0:zeros
    facs     = scanl (*) 1 pos
\end{verbatim}
Differentiation and integration enjoy the appropriately succinct definitions,
\begin{verbatim}
    deriv f = zipWith (*) pos (tail f)
    integ f = 0:zipWith (/) f pos 
\end{verbatim}
The definitions $D\ \exp=\exp;\ \exp_0=1$ and $x^*=1+xx^*$ have the following solutions in Haskell.  An {\tt x} is affixed to prevent name clashes with existing names (for example {\tt exp} in Haskell implements the function $e^x$).  One can test $x^* = \leo \exp$ by checking the first few terms of their difference.
\begin{verbatim}
    expx  :: (Eq a,Fractional a) => [a]
    expx  = 1 + integ expx
    starx :: (Eq a, Num a) => [a]
    starx = 1 : starx
    > takeW 6 (starx - e2o expx)
    [0,0,0,0,0,0]
\end{verbatim}
A rational $a/b$ is presented in Haskell as {\tt a\%b}.  The elements of {\tt expx} are rationals, and {\tt e2o} removes the factorial divisors, yielding {\tt [1\%1,1\%1, \ldots ]}.  The following defines {\tt takeW n} which is {\tt take n} preceded by the conversion of whole rationals into integers (the (.) is function composition, and {\tt properFraction} is in the Prelude).
\begin{verbatim}
    makeWhole r  = case properFraction r of
                    (n,0)         -> n
                    otherwise     -> error "not whole"
    makeAllWhole = map makeWhole 
    takeW n      = take n . makeAllWhole
\end{verbatim}
Table \ref{tab-core} contains further core sequences defined by differential equations.  All should be given the type {\tt (Eq a, Fractional a) => [a]}, like {\tt expx} above.  The core sequence $x^*$, which we defined earlier, could be defined by \verb|starx = 1+integ (starx^2)|, since $D\,x^* = (x^*)^2;\ x^*_0=1$ (but its elements would then be fractional).  For two more examples, let us calculate definitions for $\arctan$ and $\arcsin$.  
\[ D\, \conv{\tan}=((1+\tan^2)\co \conv{\tan})^{-1}=1/(1+x^2) \]
\[ D\,\conv{\sin} = (\cos \co \conv{\sin})^{-1} 
=  \left(\left(\sqrt{1-\sin^2}\right) \co \conv{\sin}\right)^{-1}
=  1/\sqrt{1-x^2} \]
The latter uses the Pythagorean identity, $\sin^2+\cos^2=1$, which follows from the defining differential equations.  Taking initial values into consideration, the solutions rendered in Haskell are immediate:
\[\begin{array}{l}
\verb|atanx   = integ (1/(1+x^2))|\\
\verb|asinx   = integ (1/(sqroot (1-x^2)))|
\end{array}\]
Here are checks of $\exp = (\sec+\tan) \co \gd$ ($\gd$ is the Guddermanian function) and $D\,\conv{\sin} =  1/\sqrt{1-x^2}$.
\begin{verbatim}
    > takeW 6 (expx - ((secx + tanx) `o` gdx))
    [0,0,0,0,0,0]
    > takeW 6 (deriv (converse sinx) - (1/(sqroot (1-x^2))))
    [0,0,0,0,0,0]
\end{verbatim}
\begin{table}
\[\begin{array}{l|l}
\verb|lgnx    = integ (1/(1+x))| &       D\, \lne = 1/(1+x);\ \lne_0=0\\
\verb|sinx    = integ cosx| &		    D \sin = \cos;\ \sin_0=0\\
\verb|cosx    = 1-integ sinx| &			D \cos = -\sin;\ \cos_0=1\\
\verb|tanx    = integ (1+tanx^2)| &		D \tan = 1+ \tan^2;\ \tan_0=0\\
\verb|secx    = 1+integ (secx * tanx)| & D \sec = \sec\tan;\ \sec_0=1 \\
\verb|coshx   = 1+integ sinhx| &		D \cosh = \sinh;\ \cosh_0=1\\
\verb|sinhx   = integ coshx| &			D \sinh = \cosh;\ \sinh_0=0\\
\verb|tanhx   = integ (1-tanhx^2)| &	D \tanh = 1-\tanh^2;\ \tanh_0=0\\
\verb|gdx     = integ (1/coshx)| & D\,\gd = 1/\cosh;\ \gd_0=0 
\end{array}\]\caption{Some core Sequences\label{tab-core}}
\end{table}

\noindent A bivariate sequence, $b(z,u)$, may be regarded as a (potentially doubly-infinite) matrix, $t=(b_{i,j})$, of coefficients of $z^iu^j$.  It is implemented as a univariate sequence, $s$, of (homogeneous) polynomials such that $s_n$ is the diagonal, $[b_{0,n},\ b_{1,n-1}, \ldots b_{n,0}]$ of $t$.  Thus, $b_{0,n}z^0u^n+b_{1,n-1}z^1u^{n-1}+\cdots+b_{n,0}z^nu^0$ is represented by $s_n=b_{0,n}x^0+b_{1,n-1}x^1+\cdots+b_{n,0}x^n$ ($z$ becomes $x$, and $u$ is redundant).  So, $[u^{n-k}z^k]b=[x^k][x^n]s$; equivalently $[u^kz^n]b=[x^n][x^{n+k}]s$. The following depicts the $t$ to $s$ map on a portion of $t$:
\[ 
\left[\begin{array}{rrrrr}
b_{0,0} & b_{0,1} & b_{0,2} &b_{0,3} & \ldots\\
b_{1,0} & b_{1,1} & b_{1,2} &b_{1,3}& \ldots\\
b_{2,0} & b_{2,1} & b_{2,2} & b_{2,3} &\ldots \\
b_{3,0} & b_{3,1} & b_{3,2} & b_{3,3} & \ldots
\end{array}\right] \mapsto \left[\begin{array}{rrrr}
b_{0,0} \\
b_{0,1} & b_{1,0} \\
b_{0,2}& b_{1,1} & b_{2,0} \\
b_{0,3} & b_{1,2} & b_{2,1}&b_{3,0}
\end{array}\right] 
\]
The ring-with-division and square root operations pertaining to $b(z,u)$ are isomorphically transferred to the diagonal representation $s$.  The representations for $u$ and $z$ are, $u=0u^0z^0+(1u^1z^0+0u^0z^1)\cong [0x^0,1x^0+0x^1]=[[0],[1,0]]$ and $z=0+(0u+1z)\cong [0,0+1x]=[[0],[0,1]]$.  Here is $(u+z)^*$ represented by $s$={\tt pascal} in Haskell,
\begin{verbatim}
    u,z, pascal :: (Eq a, Num a) => [[a]]
    u      = [[0],[1,0]]
    z      = [[0],[0,1]]
    pascal = starx `o` (u+z)
    > take 6 pascal
    [[1],[1,1],[1,2,1],[1,3,3,1],[1,4,6,4,1],[1,5,10,10,5,1]]
\end{verbatim}
This displays $[u^{n-k}z^k](u+z)^*=[x^k][x^n]s=\binom{n}{k}$.  Commonly, we want to show a portion of $t$, where $[x^k][x^n]t=b_{n,k}=[x^n][x^{n+k}]s$, or perhaps more commonly, such that $[x^k][x^n]t=n![u^kz^n]b=n![x^n][x^{n+k}]s=[x^n]\leo[x^{n+k}]s$.  For example, let
\[ b=\exp\co (z+uz)=(\exp\co z)(\exp \co uz) = \sum_n\sum_{k=0}^n \frac{z^{n-k}}{(n-k)!}\frac{u^kz^k}{k!}=\sum_n\sum_{k=0}^n \binom{n}{k}\frac{u^kz^n}{n!} \]
Then $n![u^kz^n]b=\binom{n}{k}$.  The functions {\tt unDiag} and {\tt unDiage2o} transpose the $s$-representation into the desired $t$-representation (the reverse of the $t\mapsto s$ map depicted above).  The function {\tt unDiage2o} first removes factorial divisors associated with $z$.  The functions {\tt select} and {\tt selectW} take an argument $[n_0,n_1, \ldots, n_m]$ saying what length to take from rows 0 to $m$ of $t$.  The version {\tt selectW} converts whole rationals to integers.  Bivariate counting sequences, $b_{n,k}$, are typically zero for $k>n$, and from these we select a lower triangular section.  The {\tt schroeder} sequence from section \ref{sec-survey}, item {\bf C}, is an example which is ordinary in $u$ and $z$, whilst {\tt ebinom} below is exponential in $z$.  The sequence {\tt powerSums} of polynomials for summing powers has a polynomial of order $n+1$ at position $n$, so a lower trapezium section $[2..4]$ is selected.
\begin{verbatim}
    list, pluralList :: (Eq a, Num a) => [a]
    list       = starx
    pluralList = list - x - 1
    schroeder  =  z + u*(pluralList `o` schroeder)
    > select [1..6] (unDiag schroeder)
    [[0],[1,0],[0,1,0],[0,1,2,0],[0,1,5,5,0],[0,1,9,21,14,0]]    
    ebinom =  expx `o` (z+u*z)
    > selectW [1..6] (unDiage2o ebinom)
    [[1],[1,1],[1,2,1],[1,3,3,1],[1,4,6,4,1],[1,5,10,10,5,1]]
    powerSums  = ((expx `o` (u*z))-1)/((expx `o` z)-1)
    > select [2..4] (unDiage2o powerSums)
    [[0 % 1,1 % 1],[0 % 1,(-1) % 2,1 % 2],[0 % 1,1 % 6,(-1) % 2,1 % 3]]
\end{verbatim}
A number of supporting functions are needed.  The {\tt unDiag} function expects its argument to be perfectly triangular, so {\tt padTri} is first applied to fill out the triangle with zeros if necessary.  Then {\tt transpose m} detaches the heads of the rows of {\tt m}, which make up the first column, {\tt c}, and this becomes the first row of the result.  A recursive invocation transposes the remaining sub-matrix, {\tt m'}.
\begin{verbatim}
    select s t   = zipWith take s t
    selectW s t  = zipWith takeW s t
    unDiag       :: Num a=> [[a]]->[[a]]
    unDiag       = transpose . padTri
    unDiage2o    :: Fractional a=> [[a]]->[[a]]
    unDiage2o    = unDiag . (map e2o)
    padTri t     = zipWith padRight t [1..]
    padRight r k = r++(take (k-(length r)) zeros)    
    
    transpose []   = []
    transpose m    = c:transpose m'
      where (c,m') = foldr detachHead ([],[]) m
            detachHead   ([r0]) b  = (r0:fst b,snd b)
            detachHead   (r0:r') b = (r0:fst b,r':snd b)
\end{verbatim}
There is plenty of room for adding functions to taste.  Perhaps the main difficulty is deciding on an effective naming convention.  Here are some examples.
\begin{verbatim}
    takeEBivW r = (selectW r) . unDiage2o
    takeEBiv  r = (select r)  . unDiage2o
    takeBivW  r = (selectW r) . unDiag
    takeBiv   r = (select r)  . unDiag
\end{verbatim}
Differentiation with respect to $z$ (or $u$) is performed by {\tt dz} (or {\tt du}).  Below are the definitions, plus a test based on the set partitions recurrence from section \ref{sec-survey}, item {\bf L}. Instead of using {\tt allZeros [1..6]}, the reader may wish to simply use {\tt selectW [1..6]} and view the zeros (incidentally, one would then observe that the diagonal representation is not always a perfect triangle).
\begin{verbatim}
    dz s = map deriv (tail s)
    du s = map (reverse . deriv . reverse) (tail (padTri s))
    set, nonEmptySet, emptySet :: (Eq a, Fractional a) => [a] 
    set          = expx
    emptySet     = 1
    nonEmptySet  = set - emptySet
    parts        = set `o` (u*(nonEmptySet `o` z))
  
    allEq c r    = foldr (\a b-> (a==c) && b) True r
    allZeros s t = allEq True (map (allEq 0) (select s t))
    > allZeros [1..6] ((dz parts) - (u*parts +u*du parts))
    True
\end{verbatim}

\section{Exercising the implementation}\label{sec-exercising}
As it stands, the implementation facilitates a great range of experimentation.  We will demonstrate a few concrete examples.  It will be seen that the implementation is a valuable assistant in the study of otherwise theoretical material.  

The definitions in tables \ref{tab-seq-exps} and \ref{tab-biv-seqs} transliterate into Haskell.  Here are some examples (see items {\bf B} and {\bf D}).  One has to be vigilant about when to use {\tt e2o, take, takeW, select, selectW, unDiag, unDiage2o}, \textit{etc}.  It is not necessary to give type declarations to all definitions, but it will be necessary for some.  We leave that as a trial-and-error exercise.
\begin{verbatim}
    cycle, perm :: (Eq a, Fractional a) => [a]
    perm  = starx
    lg g  = lgnx `o` (g-1)
    cycle = lg starx
    > takeW 6 (perm - (set `o` cycle))
    [0,0,0,0,0,0]
    cayleyTree            = x*(set `o` cayleyTree)
    connectedAcyclicGraph = cayleyTree - cayleyTree^2 / 2
    > takeW 8 (e2o connectedAcyclicGraph)
    [0,1,1,3,16,125,1296,16807]
     
    infix 7 |^
    (|^)      :: (Eq a, Fractional a) => [a] -> a -> [a]
    f |^ r    = expx `o` (r *| lg f)
    legendre  = (1-2*u*z+z^2) |^ [-1%2]
    > select [1..4] (unDiag legendre)
    [[1 % 1],[0 % 1,1 % 1],[(-1) % 2,0 % 1,3 % 2],
    [0 % 1,(-3) % 2,0 % 1,5 % 2]]
    hermite   = expx `o` (2*u*z-z^2)
    > select [1..4] (unDiage2o hermite)
    [[1 % 1],[0 % 1,2 % 1],[(-2) % 1,0 % 1,4 % 1],
    [0 % 1,(-12) % 1,0 % 1,8 % 1]]
\end{verbatim}
No attention has been paid to efficiency or prettiness of results -- all the computations are expected to work on small examples, resulting in small results.  Endless examples could be given related to items {\bf A-Z}.  We must choose only a few, and we aim for variety.  

The factorials are defined in the previous section using {\tt scanl}; below they are generated directly from their differential equation, by a continued fraction recurrence, and by shuffle inverse (see items {\bf E} and {\bf P}).  Shuffle product can also be defined by $f\otimes g=\leo(\loe f * \loe g))$ (but that involves rationals, even when $f$ and $g$ are integer sequences).
\begin{verbatim}
    fac = 1+x*fac + x^2*(deriv fac)
    > take 6 fac
    [1,1,2,6,24,120]
    cf_fac  = 1/(cfdenom 1)
     where cfdenom n = 1 - x*(2*n-1)-x^2*n^2*(1/(cfdenom (n+1))) 
    > takeW 6 cf_fac
    [1,1,2,6,24,120]
    infix 7 |><|    -- shuffle product
    f@(f0:f') |><| g@(g0:g') = (f0 * g0): ((f' |><| g)+(f |><| g'))
    _ |><| []                = []
    [] |><| _                = []
    shInv f@(f0:f') = (1/f0): (-f' |><| ((shInv f) |><| (shInv f)))
    > takeW 6 (shInv (1-x))
    [1,1,2,6,24,120]
\end{verbatim}
The Newton transform, $({\cal N},{\cal N}^{-1})$, is an isomorphism between the Hadamard and infiltration rings (see item {\bf K}).  It is implemented here by ({\tt h2i, i2h}), with a recursive variant, {\tt rh2i}. We also translate the definitions of $\Delta$ and $\Sigma$ directly to {\tt delta} and {\tt sigma}.  Later, we shall define another version of $\Sigma$, named {\tt prefixSums}, which produces a finite result on a finite sequence.  Let's throw into our test the ubiquitous fibonacci sequence, defined by $f_{n+2}-f_{n+1}-f_n=0;\ f_0=f_1=1$.  The first part can be re-expressed $b(E)f=0$, where $b=x^2-x-1$, with solution $f=\displaystyle\frac{(\revp{b}*[1,1])[0..1]}{\revp{b}}=1/(1-x-x^2)$ (see item {\bf S}).  For illustration, we let Haskell take the last step.  
\begin{verbatim}
    h2i s          = (1/(1+x)) |><| s
    i2h s          = (1/(1-x)) |><| s
    rh2i s@(s0:s') = s0: rh2i (delta s)
    delta s        = (tail s) - s
    sigma s        = x*starx*s
    recur          :: (Eq a, Num a) => a -> [a]
    recur a        = a:recur a
    fib            = (take 2 (rb*[1,1]))/rb where rb = reverse (x^2-x-1)
    > takeW 10 (recur 1 + sigma (delta fib))
    [1,1,2,3,5,8,13,21,34,55]
    > takeW 10  (i2h (h2i fib))
    [1,1,2,3,5,8,13,21,34,55]
\end{verbatim}
The Hadamard product, \verb+|*|+, and the infiltration product, \verb+|^|+ (and {\tt infProd}) are now introduced, and testing them is left as an exercise. 
\begin{verbatim}
    infix 7 |*|
    f |*| g = zipWith (*) f g

    infix 7 |^|
    f@(f0:f') |^| g@(g0:g') = (f0*g0): ((f'|^|g)+(f|^|g'))+(f'|^|g')
    _ |^| []                = []
    [] |^| _                = []

    infProd f g = h2i (i2h f |*| i2h g)
\end{verbatim}
Translation to and from falling factorial polynomials can be exercised as follows.  There are, of course, more efficient ways of generating the data used here ({\tt cycles, parts}), but we stick to a simple translation of the mathematical definitions.  In the first test, we use the fact (item {\bf N}) that $[0,1,5,14,30,55,\ldots]$ is representable by a polynomial of degree 3.  The second test compares falling factorials and cycle numbers.
\begin{verbatim}
    cycles        = set `o` (u* (cycle `o` z))
    fall n m      = product [n-i | i<-take m nats]
    alt           :: Num a => a->[a]
    alt r         = r:alt (-r)
    altMat m      = zipWith op (alt 1) m
                    where op sign r = zipWith (*) (alt sign) r
    monom2FacPoly = takeEBiv [1..] parts
    facMonom2Poly = altMat (takeEBiv [1..] cycles)

    toFacPoly p   = sum (zipWith (*|) p monom2FacPoly)
    fromFacPoly p = sum (zipWith (*|) p facMonom2Poly)
    squaresFacPoly= o2e (take 4 (h2i [0,1,5,14,30,55]))
    > fromFacPoly squaresFacPoly
    [0 % 1,1 % 6,1 % 2,1 % 3]
    > [fall x i | i<- [0..5]] == take 6 facMonom2Poly
    True
\end{verbatim}
The Maclaurin and Taylor expansions (item {\bf Y}) can be coded and tested:
\begin{verbatim}
    maclaurin f = o2e (map head (iterate deriv f))
    taylor f    = map o2e (zp (map (`o` u) (iterate deriv f)))
                  where zp (g0:g') = g0 + z* (zp g')

    bsinx, tsinx :: [[Rational]]
    bsinx = sinx `o` (u+z)
    tsinx = taylor sinx
    > select [1..8] bsinx == select [1..8] tsinx
    True
\end{verbatim}
Let us move beyond the {\bf A-Z} items and look at some other examples.  The Logan polynomials \cite[sect. 6.5]{GrKnPa94}, have the tangent numbers as constant terms.  Here they are, defined by a closed expression, $(\sin \co z+u\cos \co z)/(\cos \co z-u \sin\co z)$, and by an iteration.
\begin{verbatim}
    logan = (((sinx `o` z)+u*(cosx `o` z))/
            ((cosx `o` z)-u*(sinx `o` z)))
    > takeEBivW [2..5] logan
    [[0,1],[1,0,1],[0,2,0,2],[2,0,8,0,6]]
    loganPolys = iterate (\p -> (1+x^2) * deriv p) x
    > take 4 loganPolys
    [[0,1],[1,0,1],[0,2,0,2],[2,0,8,0,6]]
\end{verbatim}
The Entringer triangle, $E$ \cite{Street02,Stanley10}, has the tangent numbers on the first column (disregarding the first element), and the secant numbers on the diagonal.  Below, the triangle is generated first by a backwards and forwards (boustrophedonic) computation of partial sums, then as the diagonal (homogeneous) presentation of coefficients of $A=(\sin \co u + \cos \co u)/\cos \co (u+z))$ (named \sfn{zigzags} in table \ref{tab-biv-seqs}.  The bivariate $A$ is exponential in both $u$ and $z$, and $E_{n,k}=(n-k)!k![u^{n-k}z^k]A$ can be shown  \cite[ex. 6.75]{GrKnPa94}.  Forward partial sums are prefix sums.  Earlier, in item {\bf I}, we met $\Sigma s=xx^* s$ for computing them, and we used this above to define {\tt sigma}.  But that operator always results in an infinite sequence, even when applied to a finite one.  So here we use a different definition.
\begin{verbatim}
    prefixSums      = scanl (+) 0
    suffixSums      = reverse.prefixSums.reverse
    alternate f g a = a:alternate g f (f a)
    entringer       = alternate suffixSums prefixSums [1]
    > take 7 entringer
    [[1],[1,0],[0,1,1],[2,2,1,0],[0,2,4,5,5],[16,16,14,10,5,0],
    [0,16,32,46,56,61,61]]
    zigzags =   (((sinx `o` u) + (cosx `o` u))/(cosx `o` (u+z)))
    ue2o   :: Fractional a => [a]->[a]
    ue2o   = reverse.e2o.reverse
    > selectW [1..7] (map e2o (map ue2o zigzags))
    [[1],[1,0],[0,1,1],[2,2,1,0],[0,2,4,5,5],[16,16,14,10,5,0],
    [0,16,32,46,56,61,61]]
\end{verbatim}    
The Moessner sieve generates the sequence $M_r=[1^r,2^r,3^r,4^r, \ldots]$, given a positive integer $r$. Kozen and Silva \cite{KoSi13} cite a variety of proofs including sequence-calculational \cite{Hinze08a} and coinduction methods \cite{NiRu2011}.  Let us see how easily we can implement some of the computations in \cite{KoSi13}. There it is shown that the Moessner procedure can be described as the computation of a succession of bivariate sequences, $b_n(z,u)$, usefully represented in diagonal (homogeneous) form, $s_n$, and that $[u^0z^r]b_n=[x^r][x^r]s_n=n^r$.  The sequence of triangles begins with Pascal's triangle, $b_0=p=(u+z)^*$, represented in diagonal form by $s_0$.  Then the triangle-to-triangle step, $s_n\mapsto s_{n+1}$, is: take the row $\rho=[x^r]s_n = [\rho_0,\rho_1,\ldots, \rho_r]$, representing $h_n(z,u)=\rho_0u^rz^0+\rho_1u^{r-1}z^1+\cdots +\rho_ru^0z^r$, and compute $s_{n+1}$ representing $b_{n+1}=h_n(z,1)*p=(\rho_0z^0+\rho_1z^1+\cdots +\rho_rz^r)*p$.  Here is the computation of the first three triangles, followed by the selection of $n^r=[x^r][x^r]s_n=${\tt sn!!r!!r} from the first 5 triangles.  The function {\tt x2z} converts $\rho=[x^r]s_n=\rho_0+\rho_1x+\cdots+\rho_rx^r$ to $h_n(z,1)$ in bivariate diagonal form (there are many ways of doing this).
\begin{verbatim}
    x2z rho     = sum (zipWith (\c zn->[[c]]*zn) rho zPowers) 
                  where zPowers = (iterate (*z) 1) 
    moessnerT r = iterate (\s -> (x2z (s!!r))*pascal) pascal
    > select [5,5,5] (moessnerT 4)
    [[[1],[1,1],[1,2,1],[1,3,3,1],[1,4,6,4,1]],
    [[1],[1,5],[1,6,11],[1,7,17,15],[1,8,24,32,16]],
    [[1],[1,9],[1,10,33],[1,11,43,65],[1,12,54,108,81]]]
    nats2Power r = [sn!!r!!r | sn <- moessnerT r]
    > take 5 (nats2Power 4)
    [1,16,81,256,625]
\end{verbatim}
The function {\tt moessnerT} is easily changed to one which produces the $n$-indexed sequence $[x^r]s_n$ representing $h_n(z,1)$, and this makes way for a generalisation.  The iteration step is $[x^r]s_n \mapsto [x^r]s_{n+1}$, and the iteration starts with $[x^r]s_0=1$, representing $h_0(z,1)$.
\begin{verbatim}
    moessnerH r = iterate (\rho -> ((x2z rho)*pascal)!!r) 1
    > select [5,5,5,5] (moessnerH 4)
    [[1],[1,4,6,4,1],[1,8,24,32,16],[1,12,54,108,81]]
\end{verbatim}
Kozen and Silva generalise Moessner's theorem to encompass theorems by Long and Paasche \cite{KoSi13}.  In the generalised implementation, there are two new parameters, $h_0(z,1)$ (regarded as univariate, to be converted to bivariate by {\tt x2z}), and $d$, a sequence $[d_0,d_1,\ldots]$ of non-negative integers.  The iteration step implements the following recurrence, in which the final subscript indicates the selection of the homogeneous polynomial of degree $(\pdeg h_n)+d_n$:
\[ h_{n+1}(z,u) = [h_n(z,1)*p]_{(\pdeg h_n)+d_n} \]
Thus, rather than a simple iteration, we scan along $d$, because step $n$ requires $d_n$.  Let $c_n,\ n>0$, be the leading coefficient of $h_n(z,1)$ (which we extract using {\tt last}, since the highest-order coefficient is at the end).  The generalised theorem entails: (a) when $h_0(z,1)=1$ and $d=[r,0,0,\ldots]$ we get Moessner's result, $c_n=n^r$; (b) when $h_0(z,1)=b+(a-b)z$ and $d=[r,0,0,\ldots]$, we get Long's result, $c_n=(a+(n-1)b)n^r$; and (c) when $h_0(z,1)=1$ and $d=[d_0,d_1,\ldots]$, we get Paasche's result, $c_n=\displaystyle \prod_{i=0}^{n-1}(n-i)^{d_i}$.  Here is a rather succinct implementation.
\begin{verbatim}
    ksmlp h0 d = map last (scanl step h0 d)
     where step hn dn = ((x2z hn)*pascal)!!((length hn - 1)+dn)

    moessner r = ksmlp 1 (r:zeros)
    long a b r = ksmlp [b,a-b] (r:zeros)
    paascheFac = ksmlp 1 [1,1..]
    superFac   = ksmlp 1 [1,2..]
\end{verbatim} 
The above computations convey some Haskell by example, and demonstrate a wealth of experimentation assisted by sequence operations.  The core set of definitions are kept to a minimum, so that they are manageable in one file, and should not daunt beginners.  In keeping to this principle, we have not implemented an equivalent of formal Laurent series, so we cannot accommodate a sequence for $\cot$ (and $\csc$, and so on). However, we can define $x\cot=(x\cos)/\sin$.  Then, with reference to item {\bf X}, let $c(r)=rx(\coth \co rx)$, and test $x\coth=c(1/2)\co(2x)$, $c(1/2)\co -x=c(1/2)$ and $c(i) = x\cot$ (Gaussian rationals, introducing $i$, are in the next section):
\begin{verbatim}
    xcotx, xcothx :: (Eq a, Fractional a) => [a]
    xcotx   = (x*cosx)/sinx
    xcothx  = (x*coshx)/sinhx
    xcth r  = [r]*(x*(coshx `o` ([r]*x)))/(sinhx `o` ([r]*x))
    > take 10 xcothx == take 10 ((xcth (1%2)) `o` (2*x))
    True
    > take 10 (xcth (1%2)) == take 10 ((xcth (1%2)) `o` (-x))
    True
    > take 10 xcotx == take 10 (xcth i)
    True
\end{verbatim}
Sometimes there is simply extra work to be done to convert a sequence expression into a form acceptable by our definitions.  For example, the following expression for counting permutations by number of valleys is derived in \cite{FeSi14}:
\[ K(z,u) =1-\frac{1}{u}+\frac{\sqrt{u-1}}{u}\tan\co(z\sqrt{u-1}+\arctan\co\frac{1}{\sqrt{u-1}}) \]
This fails to compute in our implementation for three reasons (can you spot them?).  But, by using the double-angle identity for $\tan$, and $i\tanh = \tan\co ix$, it can be manipulated into the following form \cite{ElNo03}, which does compute:
\[ K(z,u) = \frac{\sqrt{1-u}}{\sqrt{1-u}-\tanh\co (z\sqrt{1-u})} \]
\begin{verbatim}
    valleys = r/(r - (tanhx `o` (z*r))) where r = sqroot (1-u)
    > takeEBivW [1..6] valleys
    [[1],[1,0],[2,0,0],[4,2,0,0],[8,16,0,0,0],[16,88,16,0,0,0]]
\end{verbatim}
The question of how to circumscribe a minimal core set of definitions that perhaps manifest a timeless quality, is a challenging one.  It seems only too easy to keep adding stuff, as the next section testifies.

\section{Building on the implementation}\label{sec-extensions}
Implementing sequence algebra is an example of mathematics-programming synergy, as found for example, in \cite{NiWi78,Lipson81,StWh86,RB:88,ODHaPa08,Wilf02,SuWi12,DoVE12}.  One should note the chronology of language use: \cite{NiWi78} uses Fortran, \cite{Lipson81} uses pseudo-Algol, \cite{StWh86} uses Pascal, \cite{RB:88} uses Standard ML, \cite{Wilf02} uses Maple, \cite{SuWi12} uses Scheme, and \cite{ODHaPa08,DoVE12} use Haskell.  Haskell is one of the most recent and ambitious in the evolution of programming languages.  The story of its development \cite{HuHuPJWa07} is an informative account of collaborative design in a scientific context.  It clearly reveals the tensions between the pursuit of elegant tried-and-tested universal concepts, and pragmatically-motivated more complex and speculative features.  

One has to face the fact that Haskell presents the casual newcomer with subtleties, some of which cause bafflement.  This slightly detracts from our goal, but also means that the implementation of sequence algebra is a fine benchmark test: Haskell ought to host it well for relative beginners.  There are two prominent sources of subtleties: lazy evaluation and type classes.  The former might be discovered in working with infinite matrices, for example try rewriting {\tt transpose}.  The latter is likely to cause the most frustration.  One could write an elucidation of potential ``surprises'' centred around implementing sequence algebra.  That is beyond our scope, but we draw attention to the fact that some type declarations can be omitted, and some not.  To take just one example, the final test of the previous section, {\tt take 10 xcotx == take 10 (xcth i)}, does not go through if the type declaration for {\tt xcotx} is omitted (then the system doesn't know to translate the rationals in {\tt xcotx} to Gaussian rationals for comparison).  On the other hand, we may omit an explicit type for {\tt xcot} and use the test
\begin{verbatim}
    makeReal (r:&0)  = r
    makeReal _       = error "not real"
    makeAllReal g    = map makeReal g 
    > take 10 xcotx == makeAllReal (take 10 (xcth i))
    True
\end{verbatim}
However, if the {\tt g} is omitted from the definition of {\tt makeAllReal}, then {\tt makeAllReal} is given a different type and the test fails to type-check.  Of course, such things have interesting explanations, but they are potentially off-putting for beginners.

These remarks notwithstanding, one cannot resist adding to the implementation in a myriad of ways.  Here are a few next-steps, which the author has already taken, and which are left as fruitful exercises. 
\begin{itemize}
 \item Translate \cite{Wilf02}, and elements of \cite{StWh86}, to use Haskell.
 \item Introduce Gaussian rationals as an instance of {\tt Num} and {\tt Fractional}, and test De Moivre's theorem (item {\bf A}).  Here is part of a definition and a test of Euler's identity:
\begin{verbatim}
    infix  6  :&
    data Gaussian a  = a :& a deriving (Eq, Read, Show)
    i    :: (Eq a, Num a) => Gaussian a
    i    = 0:&1
    ix   :: (Eq a, Num a) => [Gaussian a]
    ix = [0,i]

    instance  (Num a) => Num (Gaussian a)  where
      -- define negate, +, abs, signum, fromInteger
      (x:&y) * (x':&y') =(x*x'-y*y') :& (x*y'+y*x')

    instance  (Fractional a) => Fractional (Gaussian a)  where
      (x:&y) / (x':&y') =  (x*x'+y*y') / d :& (y*x'-x*y') / d
                           where d  = x'*x' + y'*y'
      fromRational a    = fromRational a :& 0
    
    > take 10 (cosx + [i]*sinx) == take 10 (expx `o` ix)
    True
\end{verbatim}
\item Introduce an instance, {\tt Shuffle a}, of class {\tt Num}, so that one can write \verb+s |><| t+ as {\tt S s * S t}, and shuffle power $\shp{s}n$ as \verb+(S s)^n+. Extend the following code, making {\tt Shuffle a} an instance of {\tt Fractional}.  The first test below illustrates shuffle power.  The second test involves the secant numbers, $s=\leo \sec$.  These can be defined (see items {\bf E} and {\bf F}) by applying $\leo$ to the differential equation for $\sec$ to give $s'=s\otimes (\leo \tan)$, and since $\sec_0=1$ we get the Haskell \verb+secNums = 1:secNums |><| tanNums+ (where {\tt tanNums=e2o tanx}).  Contrast this to the use of {\tt S}:
\begin{verbatim}
    newtype Shuffle a = S [a] deriving (Eq, Read, Show)
    unS (S s) = s

    instance (Eq a, Num a) => Num (Shuffle a) where
      negate (S s)  = S (negate s)
      (S s) + (S t) = S (s+t)
      (S s) * (S t) = S (s |><| t)
      fromInteger n = S [fromInteger n]
      abs _         = error "abs undefined on Shuffle"
      signum _      = error "signum undefined on Shuffle"

    > takeW 6 (unS ((S starx)^2))
    [1,2,4,8,16,32]
    tanNums = e20 tanx
    secNums = 1:unS (S secNums * (S tanNums))
    > takeW 10 secNums
    [1,0,1,0,5,0,61,0,1385,0] 
\end{verbatim} 
\item Introduce matrix computations.  To keep the definitions simple, use the type {\tt [[a]]} for a matrix, and presume, controversially,  that it is used responsibly, in the sense that a matrix is presented as a list of rows of agreed length.    Transpose is already defined in our implementation (written to work also for infinite matrices).  Operations to define include determinant, characteristic polynomial, adjugate, Gaussian elimination, and different methods of inversion.  Then one can test computations in proofs of the Cayley-Hamilton theorem, and experiment with bivariate Lagrange inversion (using $2\times 2$ Jacobians).
\item  Sequences of sequences can become confused with matrices, so it is instructive to define:
\begin{verbatim}
    data Matrix a    = M [[a]] | D a deriving (Eq, Read, Show)

    instance (Eq a, Num a) => Num (Matrix a) where
      negate (M m)   = M (map (map negate) m)
      negate (D r)   = D (negate r)
      (M a) + (M b)  = M ...
      ... clauses for + and *
      
      fromInteger n = D (fromInteger n)
\end{verbatim} 
The idea is that if {\tt s} is a square matrix of type {\tt [[a]]}, then we can have {\tt M s}.  Definitions of addition, {\tt M s + M t}, and multiplication, {\tt M s * M t}, can (with dereliction of duty) assume that {\tt s} and {\tt t} are square of the same dimension.  The element {\tt D r} stands for the square diagonal matrix (of any dimension) with {\tt r} along the diagonal.  The instance definitions of addition and multiplication each require four clauses (MM, DM, MD, DD), negate has two clauses (M, D):   
\item Rewrite part III of \cite{Lipson81} to use Haskell, making good use of classes and instances to reflect the algebraic structure.  At one level this can be approached as a program translation exercise, and is rewarding in demonstrating Haskell to be a good host language.  At other levels it invites study of a good bit of theory (Euclidean domains, finite fields, Chinese Remainder Theorem, interpolation, homomorphic image schemes, Fast Fourier Transform, and Newton's algorithm applied to power series).
\end{itemize}
Further to these tried-and-tested steps, there is, of course, unlimited scope for add-ons.  Related software can be found in the Hackage repository of the Haskell website (www.Haskell.org). 
\section{Concluding remarks}
It is clear that sequence algebra serves calculus: many sequence identities foretell relationships between analytic functions; it serves combinatorics: many counting sequences for discrete structures can be derived by sequence algebra; and it serves computation: it expresses the behaviour of certain kinds of automata; it leads to interpolation methods and summation formulae, and supports program calculation.  The theory could hardly be more foundational, and constructing an implementation from scratch emphasises its concreteness, and has the potential to reinforce understanding.  

We have exercised the implementation on examples from \cite{GrKnPa94,FlSe09}, demonstrating that it makes a valuable companion to those texts.  It could be applied to other texts, for example \cite{Comtet,GoJa83,BeLaLe98,Stanley11,Stanley99}.  It can also serve as a centre-piece in a course on functional programming in mathematics.  And, indeed, the experience of typing up and experimenting with the code, confronts one with intriguing issues in programming language design.  There is zero-testing on sequence elements, which could be used to open a discussion on computability.

The on-line encyclopaedia of integer sequences \cite{OEIS} has hundreds of thousands of sequences.  The sequences we have mentioned can be found using the OEIS search facility.  It will be noticed that many of the sequences are accompanied by generating code written in various languages, including Haskell.  One may like to investigate how many OEIS entries can be expressed in the ``language'' of tables \ref{tab-seq-exps} and \ref{tab-biv-seqs}.  A Haskell interface to the OEIS is reported in \cite{Winter13}.

Needless to say, to elaborate the topic more fully, with proof details and examples, one needs a book-sized exposition (draft portions of a book may be requested from the author).  Beyond that, the obvious question is how to make a seamless progression.  A few programming-oriented suggestions are in section \ref{sec-extensions}.  On the theory side, we must acknowledge that sequence algebra is so low in the mathematical hierarchy, that it doesn't determine a narrow range of follow-up topics.  Nevertheless, we mention a few.  One is the classification of sequences, taking a lead from \cite{HaKuRu17} and \cite[Ch. 6]{Stanley97}.  Related to this is the computer algebra work done under the heading ``generating functions'' or ``holonomic functions'' \cite{HeRu11,Kauers13}.  It remains to construct a bridge from the elementary level of the present paper to the use of a computer algebra package.  

Established results on differential equations, including computer-algebraic, may be revisited with an eye to drawing out those which become particularly accessible when specialised to sequences.  One suggestion is to bring the method of characteristics as used, for example in \cite{FeSi14}, into common parlance for sequence work.  Another is to find a smooth passage from the level of the present paper to results obtained using the language of \textit{Species}, for example those in \cite{BeLaLe98,PiSaSo12} (an introduction to Species for Haskell programmers is \cite{Yorgey10}).

Various multivariate directions beckon, including formal languages \cite{BR:88,BaHaPiRu17} and multivariate Lagrange inversion \cite{Gessel87}.  We have also arrived at the threshold of analysis but we have not crossed it, except for bringing $\pi$ into item {\bf X}. It is natural to ask whether fluency in infinite sequences, as promoted here, has any bearing on how students approach Cauchy sequences and analytic functions.  Related to this is the progression from chapter 1 to chapter 2 in \cite{Henrici74} (and chapter VII of \cite{Eil:74}).  On another tack, one may use sequence algebra to motivate abstract algebra.  For example, Eilenberg \cite[ch. XVI, sect. 10]{Eil:74} gives a proof of the Cayley Hamilton theorem using module concepts, and module concepts are used in \cite{Gatto12,GaLa16,GaSc15} -- papers whose titles echo \cite{Klarner69,Klarner76}, but which involve a quantum-leap in mathematical sophistication.  As a final remark, we note that the eponymous Haskell B. Curry, also abstracted from concrete operations on formal power series \cite{Curry51}.

\subsection*{Acknowledgements}
This work originated (some years ago) when I was an occasional visitor at the University of York.  I am greatly indebted to Colin Runciman for providing that opportunity, and to Colin, Jeremy Jacob, and Detlef Plump for encouragement.  Special thanks are due to Daniel Siemssen, Patrik Jansson, Tim Sears and Peter Thiemann for comments on work related to this paper.  (Also, if you are an anonymous JFP reviewer of an earlier related paper, then my thanks to you too!)  Tim Sears has placed a version of the Haskell code on www.GitHub.com (under TimSears/SequenceAlgebra).

\bibliographystyle{plain}
\bibliography{SeqAlg}

\end{document}